\documentclass[a4paper]{article}

\usepackage[english]{babel}
\usepackage[sort]{cite}

\usepackage{lipsum}
\usepackage{amsfonts,amsmath}
\usepackage{amsthm}
\usepackage{graphicx}
\usepackage{epstopdf}
\usepackage{algorithmic}
\usepackage{adjustbox}
\ifpdf
  \DeclareGraphicsExtensions{.eps,.pdf,.png,.jpg}
\else
  \DeclareGraphicsExtensions{.eps}
\fi

\newcommand{\magenta}{\color{deepmagenta}}

\newtheorem{theorem}{Theorem}[section]
\newtheorem{lemma}{Lemma}[section]
\newtheorem{definition}{Definition}[section]
\newtheorem{remark}{Remark}[section]
\newtheorem{proposition}{Proposition}[section]
\newtheorem{corollary}{Corollary}[section]

\title{Updating Katz centrality by counting walks}

\author{Francesca Arrigo\thanks{Department of Mathematics and Statistics, University of Strathclyde, Glasgow, UK
  (\texttt{francesca.arrigo@strath.ac.uk}).}
\and Daniele Bertaccini\thanks{Department of Mathematics, University of Rome Tor Vergata, Rome, Italy
  (\texttt{bertaccini@mat.uniroma2.it}, \texttt{filippo@mat.uniroma2.it}).}
\and Alessandro Filippo\footnotemark[3]}

\usepackage{hyperref}
\usepackage{bookmark}

\usepackage{url}

\usepackage{booktabs}
\usepackage{latexsym}
\usepackage{multirow}

\newcommand{\R}{\mathbb{R}}
\newcommand{\N}{\mathbb{N}}

\newcommand{\cS}{\mathcal{S}}
\newcommand{\cE}{\mathcal{E}}
\newcommand{\cEout}{\cE^\text{out}}
\newcommand{\cEin}{\cE^\text{in}}
\newcommand{\cN}{\mathcal{N}}
\newcommand{\cF}{\mathcal{F}}
\newcommand{\dist}{{\rm dist}}
\newcommand{\evec}{{\bf e}}
\newcommand{\cvec}{{\bf c}}
\newcommand{\avec}{{\bf a}}
\newcommand{\qvec}{{\bf q}}
\newcommand{\vvec}{{\bf v}}
\newcommand{\xvec}{{\bf x}}
\newcommand{\bone}{\mathbf{1}}
\newcommand{\bzero}{{\bf 0}}
\DeclareMathOperator*{\argmin}{\arg\!\min}

\usepackage{soul}

\usepackage{xcolor}

\usepackage{algorithm}

\usepackage{multicol}

\usepackage{pgf,tikz}
\usepackage{comment}
\usetikzlibrary{hobby}
\usepackage{tikz-cd}
\usetikzlibrary{arrows}
\usetikzlibrary{cd}
\usetikzlibrary{graphs}

\tikzstyle{vertex} = [fill,shape=circle,node distance=80pt]
\tikzstyle{edge} = [fill,opacity=.5,fill opacity=.5,line cap=round, line join=round, line width=50pt]
\tikzstyle{elabel} =  [fill,shape=circle,node distance=30pt]

\pgfdeclarelayer{background}
\pgfsetlayers{background,main}

\begin{document}

\maketitle

\begin{abstract}
We develop efficient and effective strategies for the update of Katz centralities after node and edge removal in simple graphs. We provide explicit formulas for the ``loss of walks" a network suffers when nodes/edges are removed, and use these to inform our algorithms. The theory builds on the newly introduced concept of $\cF$-avoiding first-passage walks. Further, bounds on the change of total network communicability are also derived. Extensive numerical experiments on synthetic and real-world networks complement our theoretical results.
\end{abstract}

\section{Motivation}
Katz centrality is one of the most popular 
centrality indices in network science~\cite{Katz_original_article,Estrada_CN_2011,Newman_Networks_2010}, 
{ which has successfully been applied to a wide range of problems, such as, for example, the identification of potential disease genes in protein-interaction networks \cite{katz_like_centrality_and_genes} or the optimization of artificial neural networks \cite{pruning_neural_networks_with_katz}.} 
This measure's popularity is due to its interpretability in terms of walks around a graph 
but also to the easiness of its computation~\cite{matrix_functions_networks_Benzi_Boito, noferini2024efficient}. Indeed, Katz centrality can be expressed as the solution to a (usually sparse) linear system, an extensively studied problem that is accompanied by a well-developed set of numerical linear algebra software libraries; see, e.g., \cite{bertaccini2018iterative,davis2006direct,saad2003iterative} and references therein.
{ Moreover, Katz centrality can be easily extended to the setting of dynamic networks \cite{Grindrodetall_evolving_communicability}  or specialised to only consider particular types of walks \cite{grindrod2018deformed} without affecting its easiness of computation, thus making it possibly the most widely used walk-based centrality measure in the field. 
}

In this paper, we consider the problem of approximately updating the vector of Katz centrality { and its ``personalised" counterpart \cite{on_the_limiting_Benzi_Klymko}} after the sequential removal/failure of a set of nodes or edges from a simple network.  
{ 
This is particularly relevant in real life applications. 
Indeed, in many real-world networks nodes and/or edges evolve over time~\cite{robustness_internet_2001,robustness_of_bio_networks_2002,robustness_power_grid_2004}. For example, an edge being removed could correspond to road closure within a road network, while a node being deleted could model a user leaving a social media platform. 
When centrality is of interest in these and similar situations, recomputing the vector of centrality scores may not always be feasible, especially for extremely large networks. Approximating strategies like the ones presented in this paper may be devised instead.
Additionally, our results may also be useful in the context of network design and network optimization, where simulating the removal of edges is often used as a way to address issues related to the robustness of the network or the presence communication bottlenecks \cite{attack_vulnerability_Holme_2002,attack_robustness_Iyer_2013,make_it_or_break_it_2014,Updating_downdating_Arrigo_Benzi}. 
}

Several updating strategies for centrality vectors have been devised over the years \cite{updating_pagerank_Langville_Meyer_2006,updating_betweenness_2022}  and, for walk-based centralities, the problem has been extensively studied (whether directly of undirectly) from the perspective of small rank updates of matrix functions and sensitivity of their entries; see, e.g.,~\cite{sherman1950adjustment,low_rank_2018,rank_1_update_ratfun_BernsteinVanLoan,On_the_stability_PozzaTudisco,Network_sensitivity_CabreraNoscheseRaichel,Sensitivity_Schweitzer}. 
For Katz centrality in particular, a few ad-hoc strategies have been proposed. In \cite{updating_katz_Nathan_Bader} the authors proposed to update centralities via an iterative refinement method applied to the original Katz vector, while \cite{personalized_Katz_BaderNathan} developed a visiting algorithm to modify the overall count of walks up to a specific length. 

We take a different approach and devise an approximation strategy that builds on the combinatorics of walks that are lost due to node/edge removal   
{ and, more specifically, on recurrence formulas that count  \emph{first-passage walks}~\cite{Maxflow_and_counting_subgraphs_JacksonSokal}, i.e., walks between two nodes that exclusively visit the final node in their last step.}

{ From our approximations, we also derive bounds on the change in the total network communicability induced by Katz centrality~\cite{TC_Benzi_Klymko} when edges/nodes are removed from the original graph. 
These bounds, together with our formulas for the updated Katz centrality scores, could be exploited to identify the new information propagators in social networks after node/edge removal; see~\cite{zhan2017identification}. 
}

\subsection{Background and notation}
In the following, we briefly recall some standard definitions from network science and set the notation that will be used throughout. We refer the interested reader to~\cite{Estrada_CN_2011,Newman_Networks_2010} for more details.

An unweighted \emph{graph} or \textit{network} $G = (V,E)$ consists of a set of  \textit{nodes} (\textit{vertices}) $V = \{1,2,\ldots,n\}$ and a set $E\subset V\times V$ of \textit{edges} 
(\textit{links}) between them.  
If $E$ is symmetric, then $G$ is said to be \textit{undirected}; \textit{directed} otherwise. 
An undirected graph with no repeated edges nor loops, i.e., edges connecting a node to itself, is said to be \textit{simple}. 
Henceforth, unless otherwise specified, all graphs considered are assumed to be simple.

Two nodes $i,j\in V$ are said to be \textit{adjacent} if there exists an edge between them. We will denote such edge as $\{i,j\}$ or $i\sim j$. We remark that every undirected edge $i\sim j$ can equivalently be viewed as a pair of opposite directed edges $i\to j$ and $j\to i$. The \textit{degree} of node $i$, denoted by $\deg(i)$, is the number of nodes adjacent to $i$.
The \textit{adjacency matrix} $A\in\R^{n\times n}$ of $G$ is entrywise defined as
\[A_{ij} = \begin{cases}
    1 & \text{if }i\sim j \\
    0 & \text{otherwise.}
\end{cases}\] 
It follows from this definition, that $A^T=A$ whenever the graph is simple.

A \emph{walk of length $r$} from node $i$ to node $j$ is
an ordered sequence of $r+1$ nodes $i_0 = i, i_1, \ldots, i_{r}=j$
such that  $\{i_k, i_{k+1}\}\in E$ for all $k=0,1 \ldots, r-1$. A walk is said to be  \emph{closed} if $i_0=i_r$.
 A \emph{path} is a walk with no repeated vertices. 
We will say that a walk $i_0 , i_1, \ldots, i_{r}$ of length $r$ \emph{passes through} or \emph{visits} node $w$ if $i_k=w$ for some $k=0\dots,r$. Similarly, we say that a walk \emph{visits edge} $e=\{u,v\}$ if it passes through nodes $u,v$ either as $i_k=u$ and $i_{k+1} = v$ or as $i_k=v$ and $i_{k+1} = u$, for some $k=0,\ldots,r-1$.   

\begin{definition}\label{def:EEW}
Let $G$ be a connected graph and $r\in\N_{>0}$. We will call a {\rm first-passage walk (FPW) of length $r$ from $i$ to $w$}  any walk of length $r$ of the form
\[i=i_0, i_1, i_2, \ldots, i_{r-1}, i_r=w\]
such that $i_k\neq w$ for all $k=0,\ldots,r-1$.
\end{definition}
We remark that the only admissible closed FPWs are those of length $r=0$.

The following is a standard result in graph theory and can be easily proved, e.g., by induction.
\begin{lemma}\label{lemma:count_walks}
    The $(i,j)$th entry of  $A^r$, for $r\in\N$, is the number of walks of length $r$ from $i$ to $j$.
\end{lemma}

An  undirected graph is said to be \textit{connected} if, for any two nodes in $G$, there exists a walk of finite length connecting them, \textit{disconnected} otherwise. 
Throughout this paper, we will assume that all graphs are connected, unless otherwise stated. 
\begin{remark}
This simplifying assumption does not weaken the results presented in later sections. Indeed, if a graph was disconnected, then its adjacency matrix would be permutation similar to a block diagonal matrix. Our results would then apply straightforwardly to the individual diagonal blocks.        
\end{remark}

The last definition we need in this section is that of \textit{(geodesic) distance} between two nodes $i,j\in V$ in a connected graph, which is defined as 
\[\dist(i,j) =\argmin_{\displaystyle {r\in\N}}\{(A^r)_{ij}>0\}.\]

 Throughout this paper we will adopt the following notation; $\bone\in\R^n$ will represent the vector of all ones, $\bzero\in\R^n$ will represent the zero vector,  $I\in\R^{n\times n}$ will be the identity matrix, and $\evec_k$ will be its $k$th column. Further, we will denote by ${\bf a}_i=A\evec_i $ the $i$th column of $A$. 
 
\subsubsection{Walk-based centrality measures}
Centrality measures are a standard tool in network analysis. Customarily, they assign a nonnegative score to entities in a graph (e.g., nodes, edges) to quantify their importance. The larger the score, the more important the entity. In this paper, we are concerned with node centrality measures. Several concepts of centralities have been introduced over the years; see, e.g.,{ \cite[Chapter 7]{Estrada_CN_2011} and \cite[Chapter 7]{Newman_Networks_2010}}. 
One important class of measures is defined in terms of entries (or sums thereof) of the formal series 
\[
\sum_{r=0}^\infty b_r A^r = b_0I + b_1A +b_2A^2+\cdots
\]
where $b_r\geq 0$ are weights. The above expression has a nice interpretation in terms of walks taking place around the network, thanks to Lemma \ref{lemma:count_walks}, and formalizes a class of measures commonly referred to as { \textit{walk-based centrality measures}; see~\cite[Section 3]{matrix_functions_networks_Benzi_Boito} and references therein.} Depending on the choice of weights $b_r$ and on the spectral properties of $A$, the formal series may converge to a matrix function $f(A)$; see, e.g., {  \cite[Section 1.2]{higham2008functions}.} In network science, the most popular choices for the weights are $b_r =\alpha^r$, for some $\alpha>0$, and $b_r = \beta^r/(r!)$, for some $\beta>0$, which give rise to the matrix resolvent and exponential, respectively. 

In this paper, we will focus on $b_r = \alpha^r$, for some $\alpha>0$. The induced centrality measure is known as \textit{Katz centrality}. Introduced in \cite{Katz_original_article} it assigns to node $i$ the $i$th entry of the vector 
\[
\xvec:=(I + \alpha A + \alpha^2A^2 + \cdots)\bone = \sum_{r=0}^\infty \alpha^rA^r\bone
\]
where $\alpha>0$ is a damping parameter. Selecting $\alpha$ so that  $0<\alpha\rho(A)<1$, where $\rho(A)$ is the spectral radius of the adjacency matrix $A$, as we will do in the remainder of the paper, ensures convergence of the above series, and thus the vector of Katz centralities $\xvec = (x_i)$ can be seen to be the solution to the (usually sparse) linear system
\[
(I-\alpha A)\xvec = \bone.
\]

\subsection{Problem setting and outline}
The goal of this paper is to develop efficient algorithms based on the combinatorics of walks for updating the vector of Katz centrality after node/edge removal(s). 
We will assume that the removal of an edge $\{i,j\}\in E$ from a network corresponds to setting to zero the relevant entries of the original adjacency matrix. Similarly, the removal of node $i$ is modelled as the removal of all the edges adjacent to it, so that the dimension of the adjacency matrix of the new graph is the same as that of the original graph. Both removal operations correspond to rank-2 modifications of $A$;
 \[
 A - [\evec_i,\,\evec_j][\evec_j,\,\evec_i]^T \qquad \text{ and } \qquad A- [\evec_i,\,\avec_i][\avec_i,\,\evec_i]^T.
 \]

Given a set $\cS\neq \emptyset$ of vertices or edges, i.e., $\cS\subset V$ or $\cS\subset E$, we will denote by $\cvec^\cS_r\in\R^n$ the vector whose $i$th entry  $(\cvec^\cS_r)_i$ is  the number of walks of length $r>0$\footnote{When $r=0$ we will set $\cvec_0^{\cS} = \bzero$.}  originating from node $i$, ending anywhere in the network, and visiting \textit{at least one} of the elements in $\cS$, for all $i=1,\ldots,n$.
Entrywise, it holds that $A^r\bone \geq \cvec^\cS_r$. The following result follows trivially from Lemma~\ref{lemma:count_walks}.
\begin{lemma}\label{lemma:trivial}
Let $A_{\cS}$ be the adjacency matrix of the graph obtained from $G$ by removing all the elements in $\cS$, then
\begin{equation} \label{eq:trivial_walk_counting}
    \cvec^\cS_r = (A^r - A_{\cS}^r) \bone.
\end{equation}
\end{lemma}
Our algorithms will build on alternative expressions for the computation of $\cvec^\cS_r$, obtained using the concepts of first-passage walks (cf.,  Definition~\ref{def:EEW}) and of the newly introduced $\cF$-avoiding first-passage walks; see Definition~\ref{def:EEW} below.  
\begin{remark}
    The formulas derived in this paper are easily adapted to the case where the entries of  $\cvec^\cS_r\in\R^n$ are set to represent walks that end at a specific (subset of) node(s) rather than at any node. To accommodate this change, one would need to replace every instance of $\bone$ in the formulas with an indicator (one-hot encoding) vector $\mathbf{v}\in\{0,1\}^{n}$ identifying the desired targets. 
\end{remark}

\begin{remark}
    Since (entrywise) $0\leq A_\cS\leq A$, { \cite[Theorem 8.1.18]{horn2012matrix} guarantees} that $\rho(A_\cS)\leq \rho(A)$ and thus any choice of $\alpha$ that ensures $\sum_r \alpha^rA^r = (I-\alpha A)^{-1}$ will also imply that $\sum_r \alpha^r(A_\cS)^r = (I-\alpha A_\cS)^{-1}$.
\end{remark}

\bigskip
The paper is organized as follows. The next two sections will be devoted to finding  equivalent expressions for $\cvec^\cS_r$ in terms of first-passage walks. These alternative formulas will allow us to describe in closed form, and in terms of ``loss of walks", the change in Katz centrality scores and related total network communicability when nodes or edges (or both) are removed from a graph; see  Sections~\ref{sec:results_katz_centrality} and \ref{sec:tc}. 
In Section~\ref{sec:approx_procedure} we describe two algorithms that update Katz centrality after node/edge removal in using $O(m)$ operations, where $m$ is the number of undirected edges in the original graph. Numerical experiments are described in Section~\ref{sec:numerical_experiment}. Conclusive remarks and future work are discussed in Section~\ref{sec:conclusions}

\section{Counting walks through a set to edges \texorpdfstring{$\cE$}{Lg}}\label{sec:edges}
We begin by addressing the problem of counting walks that are forced to pass through a specific edge ($\cS= \cE = \{e\}\subset E$). This result will then be generalized to the case where  walks have to visit at least one element in a larger set of edges, and hence, in the next section, to the problem of counting walks that visit at least one node in a set. The key simple observation that we will exploit is that all walks that visit a given edge $e$ can be uniquely split into three parts; an initial walk that does not visit $e$, a step visiting edge $e$ (in either direction), and finally another walk that roams freely around the graph. 

\begin{proposition} \label{prop:C_k_arcs}
Let $e = \{u,v\}\in E$,  $\cE = \{e\}$, and $r>0$. It holds that 
\begin{equation}\label{eq:c_1edge}
\cvec_r^{\{e\}} = \left(
\sum_{k=0}^{r-1}(A_{\cE})^{k}\left(\evec_u\evec_v^T + \evec_v\evec_u^T\right) A^{r-k-1} 
\right)\bone,
\end{equation}
where 
$A_{\cE} = A -
\mathbf{e}_u{\mathbf{e}_v}^T - \mathbf{e}_v{\mathbf{e}_u}^T$ is the
adjacency matrix of the graph obtained after the removal of $e$ from $G$.
\end{proposition}
Before proving the result, we want to remark two things. Firstly, for any given $i\in V$, $\evec_i^T(A_{\cE})^{k}\evec_u = 0$ (resp., $\evec_i^T(A_{\cE})^{k}\evec_v=0$) for all $k<\ell_u:=\dist(i,u)$ (resp., $k<\ell_v:=\dist(i,v)$), so that the above can be rewritten entrywise as 
\begin{equation}\label{eq:c_1edge_entry}
(\cvec_r^{\{e\}})_i = \evec_i^T\left(
\sum_{k=\ell_u}^{r-1}(A_{\cE})^{k}\evec_u\evec_v^TA^{r-k-1} + 
\sum_{k=\ell_v}^{r-1}(A_{\cE})^{k}\evec_v\evec_u^TA^{r-k-1}
\right)\bone.
\end{equation}
Secondly, if $r\leq \dist(i,e) := \min\{\ell_u,\ell_v\}$, then $r-1< \dist(i,e) $ and thus both summations in the right-hand side of \eqref{eq:c_1edge_entry} are zero. 
Therefore we can assume that $r\geq \dist(i,e)+1$ whenever we are working on a specific $i\in V$.

\textit{Proof of Proposition~\ref{prop:C_k_arcs}.}
We proceed entrywise. Let $i\in V$ and let $\ell_x = \dist(i,x)$ for $x=u,v$. Suppose wlog that $\ell_u\leq\ell_v$; if $r\leq\ell_u$, then there are no walks of length $r$ originating from $i$ and visiting edge $e$, therefore $(\cvec_r^\cE)_i=0$. Similarly, for all $k=0,\ldots,r-1$ it holds that $\evec_i^T(A_\cE)^k\evec_x = 0$ for $x=u,v$ and thus equality holds in \eqref{eq:c_1edge}. 
Suppose now that $r> \ell_u$ is fixed. Every walk of such length can be uniquely split into three consecutive terms;  a walk of length $0\leq k<r$ that does not visit $e$, followed by $e$ (in either direction), and then by a walk of length $r-(k+1)$. For any $r>\ell_u$, these can be counted using $\evec_i^T(A_\cE)^k\evec_x\evec_y^T(A^{r-k-1}\bone)$, where either $x=u$ and $y=v$ or $x=v$ and $y=u$, depending on which between $u$ and $v$ is visited first. Letting $k$ vary and splitting the terms depending on which between $u$ or $v$ is first visited yields the result.\qed

Our formula, which allows us to count all walks originating from a node and ending at any node after visiting a given edge at least once, can be easily generalized to the situation where we are visiting at least one element in a larger set  of edges. 
\begin{proposition} \label{prop:C_k_arcs_generalized}
Let $\cE\subset E$ be a set of edges $|\cE|\geq 1$ and assume that, for every $e\in\cE$, $e = \{u_e,v_e\}$. 
Then, 
\begin{equation*} \label{eq:C_k_arcs_generalized}
\cvec_r^\cE = \sum_{e \in \cE } \cvec_r^{\{e\}} = \sum_{e \in \cE } 
\left( 
\sum_{k=0}^{r-1} (A_{\cE})^{k}\left( \evec_{u_e}\evec_{v_e}^T + \evec_{v_e}\evec_{u_e}^T \right)A^{r-k-1}\right)\bone,
\end{equation*}
where $A_\cE = A - \sum_{e \in \cE} (\evec_{u_e}{\evec_{v_e}}^T + \evec_{v_e}{\mathbf{e}_{u_e}}^T )$ 
is the adjacency matrix of
graph obtained from $G$ after edges in $\cE$ are removed.     
\end{proposition}
The proof of this result is in essence the same as that of Proposition~\ref{prop:C_k_arcs}, and is therefore omitted. 
\begin{remark}
From the definition of $A_\cE$ 
it immediately follows that 
    \begin{equation}\label{eq:cvecE_simplest}
      \cvec_r^\cE = 
\sum_{k=0}^{r-1} ({A_\cE})^{k}A^{r-k}\bone - \sum_{k=0}^{r-1} ({A_\cE})^{k+1}A^{r-k-1}\bone.  
    \end{equation}
\end{remark}

A natural extension of Proposition~\ref{prop:C_k_arcs_generalized} is presented in the next section, where we describe formulas to count walks that visit at least one node in a set $\cN\subset V$.

\section{Counting walks through a set of nodes \texorpdfstring{$\cN$}{Lg}}\label{sec:nodes}

In this section we consider the problem of counting walks that visit at least one node in a set. This will in turn allow us to provide formulas to update the vector of Katz centrality when a graph is updated through node(s) removal; see Section~\ref{sec:results_katz_centrality}. 
For our purposes, we identify a node $w$ with the set of edges it connects to. This will allow us to adapt the results from the previous section to our current problem. 



\subsection{One node}
We begin by considering the case of one node $w$. 
Let thus $\cN = \{w\}$ and consider the corresponding set of edges  
\[\cE =  \{\{w,v\}\in E\,:\, w\in\cN\}.\]
Then, the vector $\cvec_r^{\{w\}}$ whose entries count the number of walks of length $r$ that visit $w$ at least once is defined as \[\cvec_r^{\{w\}} := \cvec_r^\cE\]
and moreover it holds
\begin{eqnarray*}
\cvec_r^{\{w\}} = \cvec_r^\cE 
&=& \sum_{\{w,v\}\in \cE }  \left( \sum_{k=0}^{r-1} (A_{\cE})^{k}(\evec_w\evec_v^T+\evec_v\evec_w^T)A^{r-k-1}\right)\bone \\
&=& \sum_{v\in V} A_{vw} \left( \sum_{k=0}^{r-1} (A_{\cE})^{k}(\evec_w\evec_v^T+\evec_v\evec_w^T)A^{r-k-1}\right)\bone \\
&=& \sum_{k=0}^{r-1} (A_{\cE})^{k}(\evec_w\avec_w^T+\avec_w\evec_w^T)A^{r-k-1} \bone \\
&=& \sum_{k=0}^{r-1}  \left( (A_{\cE})^{k}\evec_w\evec_w^TA^{r-k}+ (A_{\cE})^{k}A\evec_w\evec_w^TA^{r-k-1}\right)\bone  \\
\end{eqnarray*}
where we exploited the fact that $\sum_{v\in V} A_{vw}\evec_v = \avec_w =A\evec_w$. Now, by introducing the notation $A_\cN := A_\cE$ for $\cE = \{\{w,v\}\in E: w\in \cN\}$ and observing that 
\[(A_{\cN})^k\evec_w = \begin{cases}
    \evec_w & \text{if }k=0, \\
    \bzero & \text{if }k>0,
\end{cases}\] 
the above rewrites as 
\begin{equation}\label{eq:cN_idxtofix}
\cvec_r^{\{w\}} =  (\evec_w^TA^{r}\bone)\,\evec_w+ \sum_{k=0}^{r-1} (\evec_w^TA^{r-k-1}\bone) (A_{\cN})^kA\evec_w.    
\end{equation}
The terms $ (A_{\cN})^kA\evec_w$ appearing in~\eqref{eq:cN_idxtofix} are of particular interest to us. For any given $k=0,\ldots,r-1$, the entries of this vector count the number of all those walks of length $k+1$ that consist of a walk of length $k$ that does not visit $w$ at any point, followed by one step to reach node $w$, if such step exists.  In fact, these are all the first-passage walks (FPWs) from any node to $w$ in $G$; cf.~Definition~\ref{def:EEW}. Being able to count  FPWs efficiently, especially when nodes are sequentially removed one by one, is at the heart of the algorithms proposed in Section~\ref{sec:approx_procedure}.  

For any given walk-length ${ k}\in\N$, we will denote by $(\qvec_{k})_i$ the number of FPWs of length $k$ originating from $i$ (and ending at node $w$). 
For $i\neq w$, it holds that 
\begin{equation}\label{eq:count_FPWs1}
(\qvec_{k})_i = \begin{cases}
    0 & \text{if } 0\leq  k < \dist(i,w), \\
    (A_\cN^{k-1}A\evec_w)_i & \text{if }k\geq \dist(i,w).
\end{cases}
\end{equation}
On the other hand, when $i=w$, we have
\begin{equation*} 
(\qvec_k)_w =\begin{cases}
    1 & \text{ if } k=0, \\
    0 & \text{ otherwise.}
\end{cases}    
\end{equation*}

Using this notation, \eqref{eq:cN_idxtofix} can be rewritten as 
\[
    \cvec_r^{\{w\}} =  (\evec_w^TA^{r}\bone)\,\evec_w+ \sum_{k=0}^{r-1} (\evec_w^TA^{r-k-1}\bone) \qvec_{k+1} = (\evec_w^TA^{r}\bone)\,\evec_w+ \sum_{k=1}^{r} (\evec_w^TA^{r-k}\bone) \qvec_{k}, 
\]
and thus we have proved the following.
\begin{proposition} \label{prop:C_k}
Let $\cvec_r^{\{w\}}$ be the vector whose entries count the number of walks of length $r$ ending at any node after visiting  $\cN = \{w\}$ at least once;  then 
\begin{equation} \label{eq:C_k}
 \cvec_r^{\{w\}} =  \sum_{k=0}^{r} (\evec_w^TA^{r-k}\bone) \qvec_{k}.
\end{equation}
\end{proposition}
The expression in \eqref{eq:C_k} can be written entrywise for all $i\in V$ as 
\[
(\cvec_r^{\{w\}})_i = (\evec_w^TA^{r}\bone)\,\delta_{iw}+ \sum_{k=\ell_i}^{r} (\evec_w^TA^{r-k}\bone)(\qvec_{k})_i,
\]
where $\ell_i = \dist(i,w)$.

\subsection{More than one node}
We now want to consider the problem of counting walks that visit at least one node in a larger set of nodes $\cN\subset V$ with $|\cN| > 1$.  
We begin by observing that, differently from what was done when we removed multiple edges, in general we cannot simply set
\[
\cvec_r^\cN := \sum_{w\in\cN} \cvec_r^{\{w\}}. 
\]

\begin{figure}[t]
    \centering
    \begin{tikzpicture}[scale=.9]
        \tikzset{vertex/.style = {shape=circle,draw,minimum size=1.5em}};
        \tikzset{edge/.style = {-,thick}};
            \node[vertex] (1) at (-1.5,-1.5) {$1$};
            \node[vertex] (2) at (-1.5,1.5) {$2$};
            \node[vertex] (4) at (1.5,-1.5) {$4$};
            \node[vertex] (3) at (1.5,1.5) {$3$};
            \node[vertex] (5) at (0,0) {$5$};
        \draw[edge] (1) to[] (2);
        \draw[edge] (2) to[] (3);
        \draw[edge] (3) to[] (4);
        \draw[edge] (4) to[] (1);
        \draw[edge] (1) to[] (5);
        \draw[edge] (2) to[] (5);
        \draw[edge] (3) to[] (5);
        \draw[edge] (4) to[] (5);
    \end{tikzpicture}
    \caption{Toy network.}
    \label{fig:toy_ex}
\end{figure}
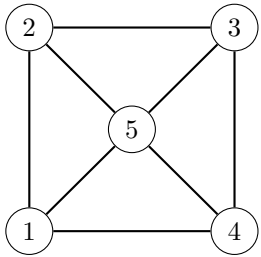

Indeed, consider the graph in Figure~\ref{fig:toy_ex} and suppose that we want our walks to visit at least one node in $\cN = \{1,2\}$. 
Using the above formula, we get, for example, $(\cvec_1^\cN)_1 = 4$; however, the number of walks of length one originating from  node 1 and visiting at least one of the nodes in $\cN$ is three ($1\to 2$, $1\to 4$, $1\to 5$). 
This mismatch arises because the nodes that we want our walks to visit are connected to each other, causing the walk $1\to 2$ to be counted twice in the right hand side of the above expression, once by the contributing term $\cvec_1^{\{1\}}$ and once by the contributing term $\cvec_1^{\{2\}}$.
This issue carries through to longer walks. Thus, a different approach should be devised that correctly counts walks when neighbouring nodes are removed. In the following, we provide two. 

The first approach uses the relationships derived in the previous section to define $\cvec_r^\cN := \cvec_r^\cE$, where, given a set of nodes $\cN = \{w_1,w_2,\ldots, w_s\}$ with $s\geq 1$, the corresponding set of edges is
\[\cE =\{\{w_i,j\}\in E:\, \text{for all } w_i\in\cN\}.\]

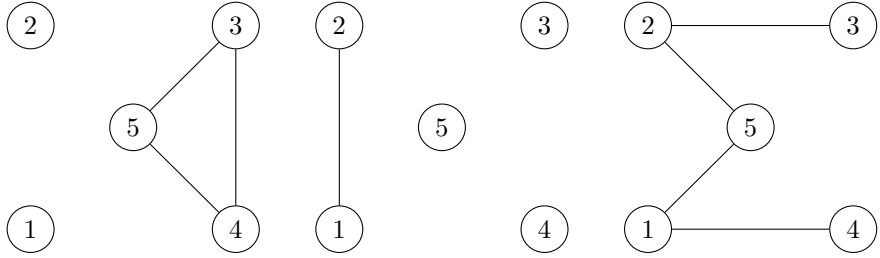
\begin{figure}[t]
    \centering
        \begin{center}
    \begin{minipage}{0.33\textwidth}
        \begin{tikzpicture}[scale=.9]
    \tikzset{vertex/.style = {shape=circle,draw,minimum size=1.5em}};
    \tikzset{edge/.style = {-,> = latex}};
        \node[vertex] (1) at (-1.5,-1.5) {$1$};
        \node[vertex] (2) at (-1.5,1.5) {$2$};
        \node[vertex] (4) at (1.5,-1.5) {$4$};
        \node[vertex] (3) at (1.5,1.5) {$3$};
        \node[vertex] (5) at (0,0) {$5$};
    \draw[edge] (3) to[](4);
    \draw[edge] (3) to[](5);
    \draw[edge] (4) to[] (5);
    \end{tikzpicture}
    \end{minipage}%
    \hfill
    \begin{minipage}{0.33\textwidth}
        \begin{tikzpicture}[scale=.9]
    \tikzset{vertex/.style = {shape=circle,draw,minimum size=1.5em}};
    \tikzset{edge/.style = {-,> = latex}};
        \node[vertex] (1) at (-1.5,-1.5) {$1$};
        \node[vertex] (2) at (-1.5,1.5) {$2$};
        \node[vertex] (4) at (1.5,-1.5) {$4$};
        \node[vertex] (3) at (1.5,1.5) {$3$};
        \node[vertex] (5) at (0,0) {$5$};
    \draw[edge] (1) to[] (2);
    \end{tikzpicture}
    \end{minipage}%
    \hfill
    \begin{minipage}{0.33\textwidth}
        \begin{tikzpicture}[scale=.9]
    \tikzset{vertex/.style = {shape=circle,draw,minimum size=1.5em}};
    \tikzset{edge/.style = {-,> = latex}};
        \node[vertex] (1) at (-1.5,-1.5) {$1$};
        \node[vertex] (2) at (-1.5,1.5) {$2$};
        \node[vertex] (4) at (1.5,-1.5) {$4$};
        \node[vertex] (3) at (1.5,1.5) {$3$};
        \node[vertex] (5) at (0,0) {$5$};
    \draw[edge] (2) to[](3);
    \draw[edge] (4) to[] (1);
    \draw[edge] (1) to[] (5);
    \draw[edge] (2) to[](5);
    \end{tikzpicture}
    \end{minipage}%
    \end{center}
    \caption{Graphs corresponding to the adjacency matrix of  $A_\cN$ (left), $A^\text{in}$ (middle) and $A^\text{out}$ (right) when $\cN = \{1,2\}$ from the network in Figure~\ref{fig:toy_ex}.}
    \label{fig:toy_ex2}
\end{figure}

Before deriving a formula for this case, we introduce some notation;
\begin{itemize}
\item $\cEin = \{\{v,w\}\in\cE:\, v,w\in\cN\}$ is the subset of $\cE$ containing edges connecting two nodes in $\cN$;
\item $\cEout = \{\{w,j\}\in\cE:\, w\in\cN\ \text{ and }j\in V\setminus \cN \}$ is the subset of $\cE$ containing edges connecting a node in $\cN$ with a node in $V\setminus\cN$, i.e., it is $\cEout = \cE\setminus\cEin$;
\item $A^\text{in}$, and $A^\text{out}$ are the adjacency matrices corresponding to the graphs $G^\text{in} = (V,\cEin)$ and $G^\text{out} = (V,\cEout)$, respectively;  
\item $A_\cN$ denotes the adjacency matrix of the graph obtained after removing the edges in $\cE$ (equivalently, nodes in $\cN$).
\end{itemize}
Considering again the graph in Figure~\ref{fig:toy_ex} we have that, for $\cN = \{1,2\}$,  the sets of edges are 
$\cE = \{\{1,2\},\{1,4\},\{1,5\},\{2,3\},\{2,5\}\}$, $\cE^\text{in} =\{\{1,2\}\}$, and $\cE^\text{out} =  \{\{1,4\},\{1,5\},\{2,3\},\{2,5\}\}$. The graphs underlying matrices $A_\cN$, $A^\text{in}$, $A^\text{out}$ are represented in Figure~\ref{fig:toy_ex2}. We point out that $ A = A_\cN + A^\text{in} + A^\text{out}$.

We are now in a position to derive an explicit formula for $\cvec_r^\cN$ for $r>0$.  By Proposition~\ref{prop:C_k_arcs_generalized} it follows 
\begin{align*}
   \cvec_r^\cN = \sum_{e\in \cE}\cvec_r^{\{e\}} &=  \sum_{e\in\cEin} \cvec_r^{\{e\}} + \sum_{e\in\cEout} \cvec_r^{\{e\}} \\
            &= \sum_{\substack{\{v,w\}\in\cE \\ v,w\in\cN}}\sum_{k=0}^{r-1}  (A_{\cN})^k\left(\evec_{v}\evec_{w}^T+ \evec_{w}\evec_{v}^T\right) A^{r-k-1}\bone \\ 
            &\qquad\qquad + \sum_{\substack{\{w,j\}\in \cE \\ w\in\cN\\ j\in V\setminus \cN}} \sum_{k=0}^{r-1}   (A_{\cN})^k\left(\evec_w\evec_j^T+ \evec_j\evec_w^T\right) A^{r-k-1} \bone\\
            &= \sum_{\substack{v,w\in\cN}}A_{vw}\left(\evec_{v}\evec_{w}^T+ \evec_{w}\evec_{v}^T\right) A^{r-1}\bone \\
            &\qquad\qquad + \sum_{\substack{w\in\cN\\ j\in V\setminus \cN}} \sum_{k=0}^{r-1}   (A_{\cN})^k\, A_{wj}\left(\evec_{w}\evec_{j}^T + \evec_{j}\evec_{w}^T\right) A^{r-k-1} \bone \\
\end{align*}
where we have used the fact that $(A_{\cN})^k\evec_w = \bzero$ for all $k>0$ and for all $w\in\cN$ and $(A_{\cN})^k\evec_i = \evec_i$ for $k=0$ and for all $i\in V$. 
Using now the following relationships
\[ \sum_{\substack{v,w\in\cN}}A_{vw}\left(\evec_{v}\evec_{w}^T+ \evec_{w}\evec_{v}^T\right) = A^\text{in} \quad \text{ and }\quad 
\sum_{\substack{w\in\cN\\ j\in V\setminus \cN}} A_{wj}\left(\evec_{j}\evec_{w}^T+ \evec_{w}\evec_{j}^T\right) = A^\text{out}\]
we derive
\begin{equation}\label{eq:cNset}
\cvec_r^\cN = 
A^\text{in} A^{r-1}\bone + 
\sum_{k=0}^{r-1}   (A_{\cN})^k A^\text{out} A^{r-k-1} \bone. 
\end{equation}
Equation~\eqref{eq:cNset} shows that we can represent $\cvec_r^\cN$ as the sum of two terms; The first counts those walks of length $r$ that start with an edge connecting two nodes in $\cE$, if such an edge exists, i.e., if $\cN$ contains two neighbouring nodes. The second term, on the other hand, includes walks that avoid nodes in $\cN$ for the first $k$ steps, than travel in one step from a node in $V\setminus\cN$ to a node in $\cN$, and then take the remaining $r-k-1$ steps around the network with no restrictions. We stress that this second term also includes walks where two neighbouring nodes in $\cN$ are consecutively visited, as long as at least one step is taken before visiting the edge connecting them. 
\begin{remark}
    When there are no neighboring nodes in $\cN$, i.e., when $\cEin = \emptyset$, it holds that $A^\text{out} = A-A_\cN$ and \eqref{eq:cNset} reduces to  \eqref{eq:cvecE_simplest}, i.e., $\cvec_r^\cN = \sum_{w\in\cN}\cvec_r^{\{w\}}$.
\end{remark}

This approach shows that it is fundamental to be able to identify the first time a node in $\cN$ is visited by a walk, especially if two or more nodes in $\cN$ are connected to each other. This leads us to the introduction of the following generalization of the concept of first-passage walks, which will be at core of our second derivation of an expression for $\cvec_r^\cN$; see Proposition~\ref{prop:C_k_generalized} below. 

\begin{definition}\label{def:GEEW}
Let $G = (V,E)$ be a graph, $i,w\in V$, $r\in \N$, and $\cF\subset V\setminus\{w\}$. An \emph{$\cF$-avoiding first-passage walk} from $i$ to $w$ of length $r$ is a FPW from $i$ to $w$ of length $r$ 
\[i=i_0, i_1, i_2, \ldots, i_{r-1}, i_r=w\]
such that $i_k\not\in\cF$ for all $k=0,\ldots,r-1$. 
\end{definition}

The newly introduced $\cF$-avoiding FPWs will allow us to count walks that visit all nodes in $\cN$ while at the same time making sure that walks that visit two neighbouring nodes in $\cN$ successively are not counted twice. To achieve this, we will use in our count the $\cF_w$-avoiding FPWs, where $\cF_w := \cN\setminus \{w\}$, as we let $w\in\cN$ vary. 
To this end, we introduce the following notation. We define the vector $\qvec_r^\cF\in\R^n$ as the vector whose $i$th entry is the number of $\cF$-avoiding FPWs of length $r\in\N$ from $i$ to $w\notin \cF$.  
Clearly, there are no closed $\cF$-avoiding FPWs  unless $k=0$, in which case $(\qvec_0^\cF)_i = \delta_{iw}$.
Moreover, if $0<k<\dist(i,w)$ it immediately follows that $(\qvec_k^\cF)_{i} = 0$. 

Thus, when $i\neq w$, it holds that 
\begin{equation}\label{eq:count_FavoidingFPWs1}
(\qvec_{k}^\cF)_i = \begin{cases}
    0 & \text{if } 0\leq k < \dist(i,w), \\
    (A_\cN^{k-1}A_{\cF}\evec_w)_i & \text{if }k\geq \dist(i,w) ,
\end{cases}
\end{equation}
where $A_\cN$ and $A_{\cF}$ are the adjacency matrices of the graphs obtained after removal of nodes in $\cN = \cF \cup \{w\}$ and in $\cF$, respectively. We note in passing that this definition implies $(\qvec_{k}^\cF)_i=0$ for all $i\in \cF$ and for all $k\geq 0$.
On the other hand, when $i=w$, we have 
\[
(\qvec_k^\cF)_w = \begin{cases}
    1  & \text{if } k = 0, \\
    0 & \text{otherwise}.
\end{cases}
\]

\begin{proposition} \label{prop:C_k_generalized}
In the above notation, the number of walks of length $r>0$ from $i\in V$ to any other node and passing through at least one of the nodes in $\cN$ is
\begin{equation}\label{eq:C_k_generalized}
(\cvec_r^\cN)_i =  
\begin{cases}
\sum_{w \in \cN} \sum_{k=1}^{r} (\qvec_{k}^{\cF_w})_{i} [A^{r-k}\mathbf{1}]_w & \text{if } i\not\in\cN \\
(A^r\bone)_i & \text{otherwise,}    
\end{cases}
\end{equation}
where $(\qvec_k^{\cF_w})_i$ is the number of $\cF_w$-avoiding FPWs of length $r$ from $i$ to $w$, and $\cF_w = \cN\setminus\{w\}$.
\end{proposition}

\textit{Proof. }
The case $i\in\cN$ holds trivially, as any walk originating from $i$ is a walk that originates from $i\in\cN$ and visits at least one node in $\cN$.

Suppose now that $i\not\in\cN$, so that $\evec_i^TA^\text{in} = \bzero^T$. Using the fact that  
\[A^\text{out} = \sum_{\substack{w\in\cN\\ j\in V\setminus \cN}} A_{jw}(\evec_j\evec_w^T + \evec_w\evec_j^T),\] 
and observing that $(A_\cN)^k\evec_w = \bzero$ for all $k>0$, \eqref{eq:cNset} becomes 
\[(\cvec_r^\cN)_i = \evec_i^T \left(\sum_{k=0}^{r-1}   (A_{\cN})^k A^\text{out} A^{r-k-1} \bone\right) =  \sum_{k=0}^{r-1} \sum_{\substack{w\in\cN\\ j\in V\setminus \cN}} \evec_i^T (A_{\cN})^k \evec_{j}\, A_{jw} \, \evec_{w}^T  A^{r-k-1} \bone .\]
The walks counted in this summation start at node $i$, take $k\geq \dist(i,j) \geq0$ steps in nodes not in $\cN$ and end in a neighbor $j\in V\setminus\cN$ of a node $w\in\cN$. Using $A^\text{out}$, these walks take one step to reach  $w\in\cN$. They finally take $r-k-1\geq 0$ steps freely around the whole network. 
These can also be counted by considering all walks that visit one of the nodes in $\cN$ for the first time after $k+1$ steps, followed by $r-k-1$ steps around the whole network, i.e., the $\cF_w$-avoiding FPW of length $k+1$ from $i$ to a node $w\in\cN$ where $\cF_w=\cN\setminus \{w\}$ , as $w$ varies, followed by $r-k-1$ steps around the graph. A change of variable then yields the conclusion.
\qed
\begin{remark}
    When $\cF=\emptyset$, Definition~\ref{def:GEEW} becomes Definition~\ref{def:EEW}, the relationships in \eqref{eq:count_FavoidingFPWs1} boil down to those in \eqref{eq:count_FPWs1} with  $\cN = \{w\}$,  and Proposition~\ref{prop:C_k_generalized} reduces to Proposition~\ref{prop:C_k}.
\end{remark}

 One of the main reasons to rewrite the count in \eqref{eq:cNset} to become \eqref{eq:C_k_generalized} is the following result, which gives us an efficient way to compute the number of $\cF$-avoiding FPWs of any given length $r>0$ and for any given set $\cF$. 

\begin{proposition} \label{prop:count_FavoidingFPWs2}
Let $\cN=\cF\cup\{w\}$, where $w\in V$ and $\cF \subset V\setminus\{w\}$. Then, in the above notation,  
\[
\qvec_{k+1}^{\cF}  =  A\qvec_k^{\cF} - \sum_{j\in\cN} (\evec_j^TA\qvec_k^{\cF}) \evec_j 
, \qquad \text{for all }\, k=0,1,2,\ldots.
\]
\end{proposition}
\textit{Proof. }
Let $A_{\cN}= A - \sum_{j\in\cN}(\evec_j\avec_j^T + \avec_j\evec_j^T)$ be the adjacency matrix of the graph obtained from $G$ after the removal of all the nodes in $\cN=\cF\cup\{w\}$, then from \eqref{eq:count_FavoidingFPWs1}, for $k=1,2,\ldots$, we get
    \begin{align*}
    \qvec_{k+1}^{\cF} &= A_\cN\qvec_k^{\cF} \\
    &= A\qvec_k^{\cF} - \sum_{j\in\cN}\left(\evec_j\evec_j^TA\qvec_k^{\cF}  + \evec_j^T\qvec_k^{\cF}\avec_j \right)\\
    &= A\qvec_k^{\cF} - \sum_{j\in\cN} (\evec_j^TA\qvec_k^{\cF}) \evec_j- \sum_{j\in\cN}(\qvec_k^{\cF})_j\avec_j 
\end{align*}   
The fact that $(\qvec_k^{\cF})_j=0$ for all $j\in \cN$ yields the conclusion for $k=1,2,\ldots$. When $k=0$, a direct computation shows that 
\[
A\qvec_0^{\cF} - \sum_{j\in\cF}(\evec_j^T\avec_w)\evec_j = \avec_w - \sum_{j\in\cF}(\evec_j^T\avec_w)\evec_j,
\]
which is exactly the vector of walks of length one starting at any node, ending ad $w$, and avoiding nodes in $\cF$, i.e., $\qvec_1^{\cF}$. This concludes the proof.
\qed

\begin{corollary} \label{corol:count_FPWs2}
When $\cF =\emptyset$, for all $k=0,1,2,\ldots$ it holds
\[
    \qvec_{k+1} = (A\qvec_k) - (\evec_w^TA\qvec_k)\evec_w.
\]  
\end{corollary}

Since $\qvec_0^{\cF} = \evec_w$ for any given $w$ and any given subset $\cF\subset V\setminus\{w\}$,  Proposition~\ref{prop:count_FavoidingFPWs2} (resp., Corollary~\ref{corol:count_FPWs2}) implies that it is possible to compute $\qvec_{k+1}^{\cF}$ (resp., $\qvec_{k+1}$) for $k=0,1,2,\dots,r$, with $r$ matrix vector products involving the adjacency matrix $A$.

In the following section we show that the results obtained so far can be used to fully characterize the variation in Katz centrality for nodes (and induced network total communicability) in networks when edges (resp., nodes) are removed from the graph. This in turn will allow us to provide bounds on the decay of total communicability and thus develop numerical strategies for their update.

\section{Updating Katz centrality via walk count}\label{sec:results_katz_centrality}


We will show that the difference in Katz centrality after the removal of a set  $\cS$ of edges or nodes can be described as a weighted sum of all walks passing through such edges or nodes, hence generalizing \cite[Proposition 5.1.b]{ranking_PLOS_Bertaccini_Filippo} to the case in which $\cS=\cE\subset E$ or $\cS=\cN\subset V$. 
In the following, 
we will denote by $\xvec^\cS = (x_i^\cS)$ the vector of Katz centralities computed on the graph $G$ after removal of the elements in $\cS$.  
We begin by proving the following simple result. 
\begin{theorem}\label{thm:old_result_generalized}
In the above notation, given  a set $\cN\subset V$  of nodes and a set $\cE\subset E$ of edges, it holds that
\[\xvec- \xvec^{\cE} =
\sum_{r=0}^{\infty}\alpha^r \cvec_{r}^\cE\qquad \text{ and }\qquad 
\xvec - \xvec^{\cN} =\sum_{r=0}^{\infty}\alpha^r \cvec_{r}^\cN,
\]
where $\cvec_0^\cS = \bzero$ and $(\cvec_{r}^\cS)_i$ is the number of walks of length $r>0$ starting form node $i$ and ending at any node after visiting at least one element in $\cS = \cN, \cE$.
\end{theorem}

\textit{Proof. }
Let $\cS$ denote either  $\cN\subset V$ or $\cE\subset E$, and let $A_{\cS}$ denote the adjacency matrix of the graph obtained from $G$ by removing all the elements in $\cS$, then
\[
\xvec- \xvec^{\cS} = \sum_{r=0}^{\infty}\alpha^r A^r \bone - \sum_{r=0}^{\infty}\alpha^r A_{\cS}^r \bone = \sum_{r=0}^{\infty}\alpha^r (A^r - A_{\cS}^r) \bone. 
\]
Applying \eqref{eq:trivial_walk_counting} yields the result.
\qed

We observe that, for any given node $i\in V$, the expressions in Theorem~\ref{thm:old_result_generalized} can be rewritten as follows 
\[x_i- x_i^{\cE} =
\sum_{r=\ell_{\cE}+1}^{\infty}\alpha^r (\cvec_{r}^\cE)_i\qquad \text{ and }\qquad 
x_i- x_i^{\cN} =\sum_{r=\ell_{\cN}}^{\infty}\alpha^r (\cvec_{r}^\cN)_i,
\]
where $\ell_\cS = \min_{s\in \cS}\{\dist(i,s)\}$, with ${\rm dist}(i,e) = \min\{\dist(i,u),\dist(i,v)\}$ for $e=\{u,v\}$. This is because there are no walks of length $k<\ell_\cN$ from $i$ to any of the elements in $\cN$ and there are no walks of length $k<\ell_\cE+1$ from $i$ that visit at least one element in $\cE$; indeed, it takes $\ell_\cE$ steps for a walk to reach the closest of the nodes in the closest edge $e\in\cE$. After this, it takes one more step to traverse $e$.

In Sections \ref{sec:edges} and \ref{sec:nodes} we derived explicit expressions for $\cvec_r^\cS$ depending on whether $\cS$ was a set of nodes or edges. We will use these in the following in order to provide formulas to update the vector of Katz centralities when edges or nodes are removed from a network. 

\subsection{Removing edges}
We first consider the case in which $\cS = \cE\subset E$ is a set of edges. 

\begin{corollary}\label{cor:update_katz_generalized_arc} 
The change in the centrality score of nodes after the removal of a set of edges $\cE$ is 
\begin{equation} \label{eq:update_katz_arc_generalized}
\xvec - \xvec^{\cE} = \sum_{\{u,v\}\in \cE} \left[
\left(\sum_{k=0}^{\infty} \alpha^{k+1}(A_{\cE})^{k}\evec_{v} \right) x_{u} 
+ \left(\sum_{k=0}^{\infty}\alpha^{k+1} (A_{\cE})^{k}\evec_{u} \right) x_{v} \right],   
\end{equation}
where $A_{\cE}$ is the adjacency matrix of
$G$ after the removal of edges in $\cE$. Equivalently,
\[\xvec-\xvec^\cE = \alpha (I-\alpha A_\cE)^{-1}\sum_{\{u,v\}\in\cE}(x_u\evec_v+x_v\evec_u).\]
\end{corollary}

\textit{Proof.}
By plugging \eqref{eq:cvecE_simplest} into   Theorem~\ref{thm:old_result_generalized}, changing the order of summation, and performing a change of variables, we get  
\begin{align*}
\xvec - \xvec^\cE &= \sum_{r=0}^\infty \alpha^r \left(\sum_{k=0}^{r-1} ({A_\cE})^{k}A^{r-k}\bone - \sum_{k=0}^{r-1} ({A_\cE})^{k+1}A^{r-k-1}\bone\right) \\
&=\sum_{k=0}^\infty \sum_{r=k+1}^\infty \alpha^r (A_{\cE})^{k}(A-A_\cE)A^{r-k-1}\bone \\
&= \sum_{k=0}^\infty\sum_{r=0}^\infty \alpha^{r+k+1} (A_{\cE})^{k} \left(\sum_{\{u,v\}\in \cE}(\evec_{u}\evec_{v}^T + \evec_{v}\evec_{u}^T) \right)A^{r}\bone.    
\end{align*}
Recalling that $x_i = \evec_i^T(I-\alpha A)^{-1}\bone$, a simple algebraic manipulation yields the conclusion.\qed
\begin{corollary}\label{corollary:update_katz_arc}
When $\cE = \{e\}\subset E$ with  $e=\{u,v\}$, for all $i\in V$ it holds
\begin{equation*} \label{eq:update_katz_arc}
x_i-x_i^{\{e\}} =
\left(\sum_{k=\ell_v}^{\infty}\alpha^{k+1}{ ({A_{\cE}}^{k})_{iv}}\right) x_u +
\left(\sum_{k=\ell_u}^{\infty}\alpha^{k+1}{ ({A_{\cE}}^{k})_{iu}}\right) x_v,
\end{equation*}
where $\ell_w=\dist(i,w)$ for $w=u,v$, and 
$A_{\cE} = A -\evec_u\evec_v^T - \evec_v\evec_u^T$ is the
adjacency matrix of graph $G$ after the removal of edge $e$. Equivalently,
\[\xvec-\xvec^{\{e\}} = \alpha (I-\alpha A_\cE)^{-1}(x_u\evec_v+x_v\evec_u).\]
\end{corollary}

Corollary~\ref{corollary:update_katz_arc} allows us to provide an upper bound on the relative decrease of Katz centrality of nodes $u$ and $v$ when edge $e=\{u,v\}$ is removed from the graph. This is described in the following result.

\begin{proposition}\label{prop:bound_on_ratio_arc}
Let $\cE = \{e\}\subset E$ with $e = \{u,v\}$. Then
\[
\frac{x_u^{\{e\}}}{x_u} \leq 1 - \alpha \frac{x_v}{x_u},\quad
\text{ and } \quad
\frac{x_v^{\{e\}}}{x_v} \leq 1 - \alpha \frac{x_u}{x_v}.
\]
Moreover, if $u$ is a leaf in $G$, i.e., if $\deg(u)=1$\footnote{Recall that $G$ is connected and thus at most one of the two nodes can be a leaf, unless $G$ is the path with two nodes.}, then 
\[
\frac{x_u^{\{e\}}}{x_u} = 1 - \alpha \frac{x_v}{x_u}.
\]
\end{proposition}
\textit{Proof. }
    We prove the first inequality. The second follows in a similar way. From Corollary~\ref{corollary:update_katz_arc}, since $\ell_u = \dist(u,u) = 0$ and $\ell_v = \dist(u,v) = 1$, and using the fact that $A_\cE\geq 0$ one gets
    \[
        x_u - x_u^{\{e\}} = x_u\sum_{r=1}^\infty \alpha^{r+1}(A_\cE^r)_{uv} + x_v\sum_{r=0}^\infty \alpha^{r+1}(A_\cE^r)_{uu} \geq  x_v\alpha (A_\cE^0)_{uu} = \alpha x_v.
    \]
    The conclusion thus follows from trivial algebraic manipulations, since $x_u\geq 1$.
    Furthermore, if $u$ is a leaf connected to node $v$ via the edge that is being removed, then $(A_\cE^r)_{uv} = (A_\cE^r)_{uu} = 0$ for all $r\geq 1$. When $r=0$, it holds $(A_\cE^0)_{uu} = 1$ and thus the conclusion. 
\qed

\subsection{Removing nodes}
We consider now the situation where a set of nodes $\cN$ is removed from the graph. Similarly to Section~\ref{sec:nodes}, the removal of one or more vertices will be modelled as the removal of the edges incident to them, resulting in the nodes in $\cN$ becoming isolated nodes in the new graph.  

\begin{corollary}\label{cor:update_katz_generalized} 
Let $\cN\subset V$ be a set of nodes in $G$, then 
\begin{equation} \label{eq:update_katz_generalized}
x_i - x_i^{\cN} = 
\begin{cases}
\sum_{w\in \cN}\,
\left( x_w  \sum_{r=\ell_w}^{\infty}\alpha^r
(\qvec_r^{\cF_w})_{i}\right), & \text{ if }i\not\in\cN \\
x_i-1 & \text{otherwise}
\end{cases}
\end{equation}
where 
$(\qvec_r^{\cF_w})_{i}$ is the number of $\cF_w$-avoiding FPWs from $i$ to
$w$, with $\cF_w := \cN\setminus\{w\}$,
and $\ell_w=\dist(i,w)$.
\end{corollary}
\textit{Proof. }
The result follows, with some algebraic manipulation, from Theorem~\ref{thm:old_result_generalized} and Proposition~\ref{prop:C_k_generalized}, using the fact that $(\qvec_r^{\cF_w})_{i}=0$ for all $r<\ell_w=\dist(i,w)$ and the convention, set at the beginning of the paper, that $\cvec_0^\cN = \bzero$.
\qed

\begin{corollary}
\label{corollary:update_katz}
In the hypothesis of Corollary~\ref{cor:update_katz_generalized}, if $\cN = \{w\}$ then 
\begin{equation} \label{eq:update_katz}
x_i - x_i^{\{w\}} = \begin{cases}
    x_w\sum_{r=\ell}^{\infty}\alpha^{r}(\qvec_{r})_{i},  & \text{ if }i\neq w \\
    x_i-1 & \text{ otherwise,}
\end{cases}
\end{equation}
where $\ell = \dist(i,w)$.
\end{corollary}

Corollary~\ref{corollary:update_katz}  { provides an estimate} on the decrease rate of Katz centrality of node $i \sim w$ after the removal of $w$ from the graph.allows to provide an estimate on the 
decrease rate of Katz centrality of node $i \sim w$ after the removal of $w$ from the graph. Indeed, the following result holds. 
\begin{proposition} \label{prop:bound_on_ratio}
If $i\sim w$ then
\begin{equation*} \label{eq:bound_on_ratio}
\frac{x_i^{\{w\}}}{x_i} \leq 1 - \alpha \frac{x_w}{x_i}.
\end{equation*}
\end{proposition}
\textit{Proof. }
     Using \eqref{eq:update_katz} and the fact that $\dist(i,w) = 1$ and $(\qvec_1)_{i}=1$ since the graph is simple, one gets
    \[x_i - x_i^{\{w\}} = x_w\sum_{r=1}^{\infty}\alpha^{r}(\qvec_{r})_{i} \geq \alpha x_w. \]
    Dividing both sides by $x_i>1$ and rearranging yields the conclusion.
\qed

\section{Bounds on total network communicability}\label{sec:tc}
In the previous section we described closed formulas that allow to compute the change in Katz centrality of nodes when edges or nodes are removed from a given simple connected graph.

In this section we use the results from Corollary~\ref{corollary:update_katz_arc} and Corollary~\ref{corollary:update_katz} to find an upper bound on
the change of the \emph{total network communicability}
\cite{TC_Benzi_Klymko}, here defined as 
\[TC = \frac{\bone^T\xvec}{n},\] after the removal of a node or an edge
from a graph.
Throughout this section, we will denote by $TC^\cS$ the total communicability of the graph obtained from $G$ after removal of elements in set $\cS$.

\begin{theorem}
\label{thm:bounds_total_comm}
Given $\cN = \{w\}\subset V$ it holds that  
\begin{equation} \label{eq:bounds_total_comm}
    TC - TC^{\{w\}} \leq \frac{1}{n} (x_w^2(1-\alpha^2 {\rm deg}(w)) - 1),
\end{equation}
where ${\rm deg}(w)$ is the degree of node $w$.
\end{theorem}
\textit{Proof. }
Let us consider the set of edges $\cE =  \{\{w,v\}\in E\, :\,w\in\cN \}$. From Corollary~\ref{cor:update_katz_generalized_arc} and symmetry of $A_{\cE}$, it follows that 
\[
\bone^T(\xvec - \xvec^{\cE})= \alpha \bone^T (I- \alpha A_{\cE})^{-1}\sum_{v:\,v \sim w } \left(x_w \evec_v + x_v \evec_w\right) 
 =\alpha x_w\sum_{v \sim w }   x_v^{\cE} + x_w^{\cE} \alpha\sum_{v \sim w}  x_v .
\]
Since, by definition of Katz centrality, $x_w^{\cE}=1$ and 
\begin{equation}\label{eq:katz}
    x_w = 1 +\alpha \sum_{v \sim w }x_v,
\end{equation} we get
\begin{equation*}
  \bone^T(\xvec - \xvec^\cE) = \alpha x_w \sum_{v \sim w }  x_v^{\cE}  + (x_w - 1)
\end{equation*}
The fact that $\xvec^{\cE}=\xvec^{\{w\}}$, Proposition~\ref{prop:bound_on_ratio}, and \eqref{eq:katz} imply
\begin{align*}
\bone^T(\xvec - \xvec^{\{w\}}) &= x_w\left(1 + \alpha \sum_{v \sim w } x_v\left(\frac{x_v^{\{w\}}}{x_v}\right) \right) - 1 \\
&\leq x_w\left(1 + \alpha \sum_{v \sim w } x_v\left(1- \alpha \frac{x_w}{x_v}\right)  \right) - 1 \\
& = x_w\left(x_w - \alpha^2 \,x_w\, {\rm deg}(w)  \right) - 1  \\
& = x_w^2(1-\alpha^2\, {\rm deg}(w)) - 1,
\end{align*}
from which \eqref{eq:bounds_total_comm} immediately follows.
\qed

\begin{theorem}
\label{thm:bounds_total_comm_arc}
Let $e=\{u,v\}$ and consider $\cE = \{e\}$.  Then
\begin{equation}\label{eq:bounds_total_comm_arc}
   TC - TC^{\{e\}} \leq \frac{1}{n} \left( 2\alpha x_u x_v - \alpha^2 x_u^2  - \alpha^2 x_v^2\right).
\end{equation}
\end{theorem}
\textit{Proof. }
 From Corollary~\ref{corollary:update_katz_arc} and symmetry of $A_\cE$  
we obtain
\[
\bone^T(\xvec - \xvec^{\{e\}}) = \alpha\bone^T\left(I-\alpha A_{\cE}^{-1}\right)(x_u\evec_v + x_v \evec_u) 
 =  \alpha x_u x_v^{\{e\}}+  \alpha x_v x_u^{\{e\}}.
\]
Now by applying Proposition~\ref{prop:bound_on_ratio_arc} we get 
\begin{align*}
\bone^T(\xvec - \xvec^{\{e\}}) &= \alpha x_u x_v\left( \frac{x_v^{\{e\}}}{x_v} + \frac{x_u^{\{e\}}}{x_u} \right) \\
& \leq \alpha x_u x_v\left( 1- \alpha \frac{x_u}{x_v} + 1-\alpha \frac{x_v}{x_u} \right) \\
& { =} 2\alpha x_u x_v - \alpha^2x_u^2 - \alpha^2x_v^2,
\end{align*}
from which the conclusion straightforwardly follows.
\qed


\begin{remark}
Consider now the function $f_{\alpha}(x,y)=2\alpha xy - \alpha^2 x^2 -\alpha^2 y^2$ appearing in \eqref{eq:bounds_total_comm_arc} and assume that $0<\alpha<1$. Let $T>t>1$, then it can be showed that 
\[
    \argmin_{(x,y)\in [t,T]\times [t,T]} f_{\alpha}(x,y) =  
    \begin{cases}
     (t,t) & \text{ if } T<\frac{2-\alpha}{\alpha}t \\ 
    (t,T) & \text{ otherwise. }     
    \end{cases}
\]
Hence, the bounds in Theorem~\ref{thm:bounds_total_comm_arc} can be used to devise strategies for optimizing total communicability through edge removal, similar to the ones in \cite{Updating_downdating_Arrigo_Benzi}. Notably, if  $T<\frac{2-\alpha}{\alpha}t$, 
we find that 
\[
    \argmin_{(x,y)\in [t,T]\times [t,T]} f_{\alpha}(x,y) \, = \argmin_{(x,y)\in [t,T]\times [t,T]} x y \, = \, (t,t), 
\]
which suggests that, in order to
minimize the impact on total network communicability of an edge
removal in a network such that $\max_i{x_i}<\frac{2-\alpha}{\alpha}\min_i{x_i}$, one should remove the edge $e=\{u,v\}$ that minimizes the
product $x_ux_v$. This was in fact one of the heuristics proposed in \cite{Updating_downdating_Arrigo_Benzi} for
solving the \emph{downdate problem} considered there. 
\end{remark}

\section{Walk-induced update of Katz centrality} \label{sec:approx_procedure}
In this section, we propose algorithms to find an approximate value of  Katz centrality for the nodes in a simple graph $G$ after elements from a chosen set $\cS$ are removed. The results presented here build on the theory developed so far, and assume that the values of the original Katz centralities, stored in a vector $\xvec$, are available. The algorithms developed in this section are especially useful in the context where the elements of $\cS$ are \textit{sequentially} removed from $G$. This would normally require solving $|\cS|$ sparse linear systems, which can be burdensome when $n$ is large. On the other hand, the methods we propose will allow updating the centrality vectors with \textit{approximations} to the actual values, and only sporadically solving a sparse linear system. As we will show, this reduces the computational effort of finding Katz centrality when nodes or edges are sequentially removed from $G$, while at the same time achieving good approximations of the new centrality scores.

Throughout this section, approximate values will be denoted by $\widehat{\xvec}$.

\subsection{Removing nodes}
We begin with the problem of updating the vector of Katz centralities after a node is removed from $G$. The strategy we propose builds on the result presented in Corollary~\ref{corollary:update_katz} and simply truncates the infinite sum to $L$ terms, where $L$ is a user selected parameter. 
The proposed approximation, for $\cN = \{w\}$, thus is
\[
\widehat{\xvec}^{\{w\}} := \xvec - x_w\sum_{r=1}^L \alpha^r \qvec_{r},
\]
where 
$(\qvec_{r})_i$ is the number of FPWs of length $r$
from node $i$ to $w$. Moreover, since Katz centrality of isolated nodes is equal to one by definition, we will also have that $\widehat{x}_w^{\{w\}} = x_w^{\{w\}}=1$.
\begin{algorithm}[t]
\caption{\bf Updating Katz centrality after node removal}
\label{alg:update_katz_node}
\begin{algorithmic}[1]
           \REQUIRE $A\in \R^{n\times n}$ adjacency matrix, $\xvec\in \R^{n}$ vector of Katz centrality, $\alpha \in (0,1/\rho(A))$, $L_{\max}^{\cN} \in \N$, $\texttt{tol}\in \R$, $w$ node to be removed. 
            \ENSURE $\widehat{\xvec}^{\{w\}} \in \R^{n}$ approximate update, $L\in \N$ number of iterations to converge. 
            \STATE {$\widehat{\xvec}^{\{w\}}=\xvec$;}
            \STATE $L = 1$;
            \STATE $\qvec = \alpha* A*\evec_w$;
            \STATE $\widehat{\xvec}^{\{w\}} = \widehat{\xvec}^{\{w\}} - x_w *\qvec$;
            \WHILE {$x_w\|\qvec\|_{2}\,/\,\|\xvec\|_2  >  \texttt{tol}$ \AND  $L< L_{\max}^{\cN}$}
            \STATE $\qvec = \alpha*A*\qvec$; \label{line:compute_fpw}\label{line:mat_prod_node}
            \STATE $q_w = 0$;
            \STATE $\widehat{\xvec}^{\{w\}} = \widehat{\xvec}^{\{w\}} - x_w *\qvec$;\label{line:compute_x_w_hat} 
            \STATE $L=L+1$;
            \ENDWHILE
            \STATE $(\widehat{\xvec}^{\{w\}})_w = 1$
\end{algorithmic}
\end{algorithm}

A pseudocode for our proposed algorithm can be found in Algorithm~\ref{alg:update_katz_node}.

Given that $\xvec$, the vector of Katz centralities, is known, the computational effort to generate the approximant $\widehat{\xvec}^{\{w\}}$ is spent computing the vectors $\qvec_r$ for $r=1,\dots,L$. Thanks to the result presented in  
 Corollary~\ref{corol:count_FPWs2}, this operation only requires $L$ matrix-vector products involving the adjacency matrix $A$. 
Therefore, the choice of $L$ is crucial to obtain both a good
approximation and a reasonable computational cost. { The latter}, in particular, should be
smaller than the cost of computing Katz centrality from scratch. In our
experiments we select $L$ \textit{dynamically}.
After selecting the maximum walk-length that we want to allow, denoted by $L_{\max}^\cN$, the parameter $L$ is chosen as the smallest $0<L\leq L_{\max}^\cN$  such that 
\begin{equation*}\label{eq:stopping_criteria_nodes}
\alpha^L
x_w\frac{\|\qvec_L\|_{2}}{\|\xvec\|_{2}}<\texttt{tol}.
\end{equation*}
This is equivalent to checking that the (relative) norm of the difference between two subsequent approximations, obtained by truncating the series at $L-1$ and $L$, respectively, is smaller than a selected tolerance \texttt{tol}. 
\begin{algorithm}[t]
        \caption{\bf Updating Katz centrality after edge removal}
        \label{alg:update_katz_edge}
\begin{algorithmic}[1]
            \REQUIRE $A\in \R^{n
            \times n}$ adjacency matrix, $\xvec\in \R^{n}$ vector of Katz centrality, $\alpha \in (0,1/\rho(A))$, $L_{\max}^{\cE} \in \N$, $\texttt{tol}\in \R$, $e=\{u,v\}$ edge to be removed.
            \ENSURE $\widehat{\xvec}^{(e)}\in \R^{n}$ approximate update, $L \in \N$ number of iterations to convergence. 
            \STATE $\widehat{\xvec}^{\{e\}}=\xvec  - \alpha*x_v*\evec_u - \alpha*x_u*\evec_v$;
            \STATE $L=1$; 
            \STATE $\mathbf{s} = \alpha^2*(A*\evec_u -\evec_v)$;
            \STATE $\mathbf{t} = \alpha^2*(A*\evec_v - \evec_u)$;
            \STATE $\widehat{\xvec}^{\{e\}}= \widehat{\xvec}^{\{e\}} - x_v*\mathbf{s} - x_u *\mathbf{t}$;
            \WHILE {$\|x_v*\mathbf{s} + x_u *\mathbf{t}\|\,/\,\|\xvec\|_2  >  \texttt{tol}$ \AND  $L< L_{\max}^{\cE}$}
            \STATE $\mathbf{s} = \alpha *(A*\mathbf{s} - s_u *\evec_v - s_v *\evec_u)$; \label{line:matrix_vector_product_u_alg_2}
            \STATE $\mathbf{t} = \alpha *(A*\mathbf{t} - t_u *\evec_v - t_v *\evec_u)$; \label{line:matrix_vector_product_v_alg_2}
            \STATE $\widehat{\xvec}^{\{e\}}= \widehat{\xvec}^{\{e\}} - x_v*\mathbf{s} - x_u *\mathbf{t}$; 
            \STATE $L= L+1$;
            \ENDWHILE
\end{algorithmic}
\end{algorithm}

\subsection{Removing edges}
For the case where an edge is removed, similarly to what we have just done for the case of nodes, we exploit the result in Corollary~\ref{corollary:update_katz_arc} to derive an approximation to
$\xvec^{\{e\}}$; 
\begin{equation*} \label{eq:approx_katz_edges}
\widehat{\xvec}^{\{e\}} := \xvec -
x_{u}\sum_{r=0}^{L}\alpha^{r+1}(A_\cE)^r\evec_v -
x_{v}\sum_{r=0}^{L}\alpha^{r+1}(A_\cE)^r\evec_u,
\end{equation*}
where $A_\cE = A - \evec_u{\evec_v}^T -
\evec_v{\evec_u}^T$ is the adjacency matrix after the
removal of edge $e = \{u,v\}$. 

As in the previous case, 
we set a maximum walk length $L_{\max}^{\cE}$ and then select $0<L\leq L_{max}^{\cE}$ dynamically as the smallest index such that
\begin{equation*}\label{eq:stopping_criteria_edges}
\alpha^{L+1}\left( x_{u}
\frac{\|A_\cE^{L}\evec_{v}\|_{2}}{\|\xvec\|_{2}} +
x_{v} \frac{\|A_\cE^{L}\evec_{u}\|_{2}}{\|\xvec\|_{2}}
\right) < \texttt{tol}.
\end{equation*}

Algorithm~\ref{alg:update_katz_edge} contains the pseudocode for the proposed procedure.

\subsection{Computational cost} The following holds. 
\begin{proposition} \label{prop:computational_costs}
{ For a sparse network,} the cost of Algorithms~\ref{alg:update_katz_node} and \ref{alg:update_katz_edge} is $\mathcal{O}(Lm)$ operations, where $m$ is the number of undirected edges in the original graph.
\end{proposition}
\textit{Proof. }
For Algorithm~\ref{alg:update_katz_node}, at every step of the loop we perform a matrix-vector product on line \ref{line:mat_prod_node} and $n$ multiplications on lines \ref{line:compute_fpw} and \ref{line:compute_x_w_hat}. Thus, the computational complexity of Algorithm~\ref{alg:update_katz_node} can be estimated as $\mathcal{O}(LC + 2Ln)$ operations, where
$C$ is the cost of performing one matrix-vector product. 
Similarly, one can deduce that Algorithm~\ref{alg:update_katz_edge} 
requires $\mathcal{O}(2LC + 4Ln)$ operations. 
The conclusion follows by recalling the reasonable assumption that $C = \mathcal{O}(m)$ for sparse matrices with $2m$ nonzero entries.
\qed


Clearly, 
Algorithms~\ref{alg:update_katz_node} and \ref{alg:update_katz_edge} can be generalized to the case where a larger set of edges (resp., nodes) is \emph{simultaneously} removed from the graph. To do so, it is sufficient to truncate the infinite series in 
\eqref{eq:update_katz_arc_generalized} (resp., \eqref{eq:update_katz_generalized}) to $L$ terms. In this case it is easy to see that, under the same assumption of Proposition~\ref{prop:computational_costs}, the generalized algorithms would have a computational complexity in time of $O(|\cS|Lm)$, where $\cS$ is the set of elements being removed. In addition, since \eqref{eq:update_katz_arc_generalized} and \eqref{eq:update_katz_generalized} are both described as sums of $|\cS|$ independent terms, one could use a parallel implementation on, e.g., $|\cS|$ processors, to compute each term independently, hence achieving the theoretical computational complexity of $O(Lm)$ per processor.

\section{Numerical experiments}\label{sec:numerical_experiment}
In this section, we compare the performance of Algorithm~\ref{alg:update_katz_node}~and~\ref{alg:update_katz_edge} with other methods that, in contrast, simply recompute the centrality scores from scratch. In our implementation the damping factor in Katz centrality was set to $\alpha = 0.85/\rho(A)$, where $\rho(A)$ is the spectral radius of $A$ \cite{TC_Benzi_Klymko}. In this section, we only considered connected networks for simplicity. 

We performed experiments on both synthetic and real-world networks. For the synthetic ones, we generated $30$ instances of networks with \texttt{n} nodes  
using the functions \texttt{pref(n,5)} and
\texttt{erdrey(n,5n)} from the CONTEST Toolbox~\cite{CONTEST}. We refer the interested reader to the original reference for further details on parameter selection. Some simple statistics on the generated graphs, including number of nodes $n$, average number of edges $m$, average diameter ${\rm diam}$, average eccentricity $\epsilon$, and average condition number $\kappa(I-\alpha A)$ of $(I-\alpha A)$  are reported in Table~\ref{tab:syntethic_networks_info}. Recall that the diameter of a network is ${\rm diam} = \max_{ij}\dist(i,j)$, while the average eccentricity is  $\epsilon = \frac{1}{n}\sum_i \max_{j\in V}\{\dist(i,j)\}$. 

The  real-world networks in our dataset, both available at \cite{SuiteSparse}, are a road network of the Minnesota state (\texttt{minnesota}) and a
communication network between a group of Autonomous Systems (AS) in
the period between 1997-2000 (\texttt{as-735}). A summary of these networks' basic features is provided in Table~\ref{tab:real_networks_info}.

All experiments were run via MATLAB R2023b on a laptop PC equipped with an Intel i7-1260P CPU with a clock rate of 2.10 GHz, and 16 GB of RAM.

\begin{table}
 	\centering
	\label{tab:syntethic_networks_info}
	\caption{Synthetic networks data: number of nodes $n$, number of edges $m$, network diameter diam, average eccentricity $\epsilon$, and condition number $\kappa(I-\alpha A)$ of  $I - \alpha A$ for $\alpha = 0.85/\rho(A)$. All the values are averaged over $30$ runs.}
	\begin{tabular}{|l|c|c|c|c|c|}       
		\hline
		Network & $n$ & $m$& diam & $\epsilon$ & $\kappa(I-\alpha A)$ \\
		\hline
		\multirow{4}*{\texttt{erdrey(n,5n)}}  & 3200 & 16000 & 6.1 & 5.2 & 10 \\
		& 6400 & 32000 & 6.9 & 5.7 & 10 \\ 
		& 12800 & 64000 & 7.2 & 6.0 & 10\\
		& 25600 & 128000 & 8.0 & 6.5 & 10 \\
		
		\hline
		\multirow{4}*{\texttt{pref(n,5)}}  & 3200 & 15800 &  5.0 & 4.7  & 10.1\\
		& 6400& 31800 &5.7& 4.8 & 11.1\\ 
		& 12800 & 63700 & 6.0 & 4.9 & 11.6 \\
		&  25600 & 127000 & 6.0 & 5.2  & 11.7 \\
		\hline
	\end{tabular}
\end{table}

\begin{table}[]
    \centering
    \caption{Real-world network data: number of nodes $n$, number of edges $m$, network diameter diam, average eccentricity $\epsilon$, and condition number $\kappa(I-\alpha A)$, for $\alpha = 0.85/\rho(A)$.}
    \begin{tabular}{|l|c|c|c|c|c|}       
		\hline
		Network & $n$ & $m$& diam & $\epsilon$ & $\kappa(I-\alpha A)$ \\
		\hline
		\texttt{minnesota}  & 2640 & 3302 & 99 & 71.7 & 12.2 \\
        \texttt{as-735} & 6474 & 13895 & 9 & 6.7 & 11.6 \\
        \hline
    \end{tabular}
    
    \label{tab:real_networks_info}
\end{table}

\bigskip
In the first experiment, we considered the synthetic networks in our dataset and removed one element $s$, either a node or an edge, selected uniformly at random from either $V$ or $E$. After constructing the adjacency matrix $A_{\cS}$ of the new graph, we compute/approximate the updated Katz centrality $\xvec^{\cS}$ in three different ways;
\begin{enumerate}
    \item[(i)] \texttt{pcg $(\xvec^{(0)} = \vvec)$}: solve $(I-\alpha A_{\cS})\mathbf{y}=\bone$ using MATLAB built-in function \texttt{pcg} (conjugate gradient) without preconditioning, initial guess $\xvec^{(0)} = \vvec$, and with a relative tolerance of $10^{-5}$; 
     \item[(ii)] method in \cite{Truncated_Katz_Foster_2001}: truncate the series $\sum_{r=0}^{\infty}(A_{\cS})^r\bone$ with a tolerance of $10^{-4}$ on the stopping criteria and a maximum number of iteration of $100$; see \cite{Truncated_Katz_Foster_2001} for details; and
    \item[(iii)] Algorithm~\ref{alg:update_katz_node} or \ref{alg:update_katz_edge} with $L_{\max}^\cS=30$ and $\texttt{tol}=10^{-4}$.
\end{enumerate}
 
Before commenting on the experimental results, some additional discussion on the stopping criteria for the above is in order. Methods (ii) and (iii) share the same stopping criterion: they stop as soon as the difference between the results of two subsequent iterations is smaller than a relative tolerance. Hence, in order to compare the two methods, it is sufficient to choose the same tolerance $\texttt{tol}=10^{-4}$. For \texttt{pcg}, on the other hand,  the stopping criterion checks the norm of the residual, and not the difference between two iterates. 
Nevertheless, we note that the
norm of the difference between two subsequent iterations of the converging sequence
of approximations $\{\xvec^{(k)}\}$ for the solution $\xvec^*$ can be bounded
by a quantity proportional to the norm of the residual of the
underlying linear system. More precisely, we recall the following.
\begin{proposition}\label{prop:stopping-criterion} 
Let $\{\xvec^{(k)} \}$ be a sequence
converging monotonically to $\xvec^*$ for $k>N$, solution of the
nonsingular system $M\xvec={\bf b}$. Then,
\begin{equation*}\label{eq:bounding_sequences}
\frac{\|\xvec^{(k-1)}-\xvec^{(k)}\|}{\|\xvec^*\|}\leq 2 \kappa(M)
\frac{\|{\bf r}^{(k-1)}\|}{\|{\bf b}\|},
\end{equation*}
where ${\bf r}^{(k-1)} = M\xvec^{(k-1)}-{\bf b}$ is the residual at the $(k-1)$th iteration.
\end{proposition}

In view of Proposition~\ref{prop:stopping-criterion}, in our experiments we used a relative
tolerance of $\texttt{tol}_{\texttt{pcg}}:=\texttt{tol}/10=10^{-5}$ for the \texttt{pcg} algorithm in method (i), where $\texttt{tol} = 10^{-4}$ is the tolerance set for methods (ii)-(iii) and the denominator comes from the observed $\kappa(I-\alpha A)\approx 10$ for all the network instances considered here; cf. Table~\ref{tab:syntethic_networks_info}.

Lastly, we want to mention in passing that we do not compare our proposed strategy with the one described in~\cite{updating_katz_Nathan_Bader} since experiments not shown here returned results comparable with those obtained by \texttt{pcg $(\xvec^{(0)} = \bzero)$}.

\bigskip 
In Tables~\ref{tab:time_experiments_nodes} and \ref{tab:time_experiments_edges}, we report the average number of iterations and average CPU times needed to compute an approximation to $\xvec^{\cS}$ using the methods described. All values are averaged over $30$ runs with the same parameters.

\begin{table}[t]
	\centering
	\label{tab:time_experiments_nodes}
	\caption{Number of iterations and CPU timings required for approximating $\mathbf{x^{\cS}}$ with methods (i)-(iii) after the removal of one node selected uniformly at random. The chosen tolerances are \texttt{tol}$=10^{-5}$ for method (i), and \texttt{tol}$=10^{-4}$ for methods (ii)-(iii). The initial guess for \texttt{pcg} is specified in brackets. All the values are averaged over $30$ runs.}
	\begin{adjustbox}{width=1.2\textwidth,center=\textwidth}
    \begin{tabular}{|l|c|cc|cc|cc|cc|}
		\hline
		\multirow{2}*{Network}                & \multirow{2}*{n}  &   \multicolumn{2}{|c|}{\texttt{pcg}($\xvec^{(0)} = \bzero$)} & \multicolumn{2}{|c|}{\texttt{pcg}($\xvec^{(0)} = \xvec$)} &  \multicolumn{2}{|c|}{method in \cite{Truncated_Katz_Foster_2001}}  & \multicolumn{2}{|c|}{Algorithm \ref{alg:update_katz_node}} \\
		&        &  iter  &   time(s)  & iter  &   time(s) & iter       &     time(s)  &  $L$   &   time(s) \\
		\hline
		\multirow{4}*{\texttt{erdrey(n,5n)}}    & 3200  &  11.0 &    1.18$\cdot 10^{-3}$    & 7.5 &  1.08$\cdot 10^{-3}$ & 45.0 &  1.75$\cdot 10^{-3}$   & \bf 6.5 & 1.27$\cdot 10^{-3}$ \\
		& 6400  & 11.0 &   6.68$\cdot 10^{-3}$    & 7.2 & 4.29$\cdot 10^{-2}$ & 45.0  & 5.63$\cdot 10^{-3}$     & \bf 4.6  &  1.08$\cdot 10^{-3}$ \\
		& 12800 & 11.0 &   1.75$\cdot 10^{-2}$     & 7.0 &  1.28$\cdot 10^{-2}$ & 45.0   &  3.39$\cdot 10^{-2}$ & \bf 4.1       & 2.20$\cdot 10^{-3}$ \\
		& 25600 & 11.0 &   2.88$\cdot 10^{-2}$      & 6.7 & 1.99$\cdot 10^{-2}$ & 45.0  &  5.36$\cdot 10^{-2}$        & \bf 3.7 &  2.47$\cdot 10^{-3}$ \\
		\hline
		\multirow{4}*{\texttt{pref(n,5)}}       & 3200 & 11.4 &    1.11$\cdot 10^{-3}$        & 5.9 & 9.28$\cdot 10^{-4}$ & 45.0 &   1.74$\cdot 10^{-3}$        & \bf 4.8 &  8.32$\cdot 10^{-4}$ \\
		& 6400 & 11.8 &   6.64$\cdot 10^{-3}$        & 5.3 &  3.30$\cdot 10^{-3}$ & 45.0 &  4.78$\cdot 10^{-3}$      & \bf 3.8  & 1.09$\cdot 10^{-3}$ \\
		& 12800 & 12.1 &   1.75$\cdot 10^{-2}$        & 4.1 & 8.57$\cdot 10^{-3}$ & 45.0 &   3.21$\cdot 10^{-2}$       &  \bf 2.6  & 1.38$\cdot 10^{-3}$ \\
		& 25600 & 12.8 &   3.10$\cdot 10^{-2}$        & 3.2 & 1.28$\cdot 10^{-2}$ & 45.0 &  5.18$\cdot 10^{-2}$        &  \bf 2.5  &  2.42$\cdot 10^{-3}$ \\
		\hline
	\end{tabular}
    \end{adjustbox}
\end{table}

\begin{table}[t]
	\centering
	\label{tab:time_experiments_edges}
    	\caption{Number of iterations and CPU timings required for approximating $\mathbf{x^{\cS}}$ with methods (i)-(iii) after the removal of one edge selected uniformly at random. The chosen tolerances are \texttt{tol}$=10^{-5}$ for method (i), and \texttt{tol}$=10^{-4}$ for methods (ii)-(iii). The initial guess for \texttt{pcg} is specified in brackets. All the values are averaged over $30$ runs.}
        \begin{adjustbox}{width=1.2\textwidth,center=\textwidth}
	\begin{tabular}{|l|c|cc|cc|cc|cc|}
		\hline
		\multirow{2}*{Network}                & \multirow{2}*{n}  &   \multicolumn{2}{|c|}{\texttt{pcg}($\xvec^{(0)} = \bzero$)} & \multicolumn{2}{|c|}{\texttt{pcg}($\xvec^{(0)} = \xvec$)} &  \multicolumn{2}{|c|}{method in \cite{Truncated_Katz_Foster_2001}}  & \multicolumn{2}{|c|}{Algorithm \ref{alg:update_katz_edge}} \\
		&        &  iter  &   time(s)  & iter   &   time(s) & iter       &     time(s)  &  $L$   &   time(s) \\
		\hline
		\multirow{4}*{\texttt{erdrey(n,5n)}}    & 3200  & 11.0 &    1.19$\cdot 10^{-3}$     & 6.1 & 1.06$\cdot 10^{-3}$ & 45.0 &  1.71$\cdot 10^{-3}$      & \bf 2.9 & 5.36$\cdot 10^{-4}$ \\
		& 6400  & 11.0 &   5.19$\cdot 10^{-3}$     & 5.8 & 3.35$\cdot 10^{-3}$ & 45.0  &  5.49$\cdot 10^{-3}$      &\bf 2.4 &  4.56$\cdot 10^{-4}$ \\
		& 12800 & 11.0 &   1.38$\cdot 10^{-2}$   & 5.5 & 8.65$\cdot 10^{-3}$ & 45.0   &  2.50$\cdot 10^{-2}$ & \bf 2.1     & 7.75$\cdot 10^{-4}$ \\
		& 25600 & 11.0 &  2.85$\cdot 10^{-2}$      & 5.1 & 1.68$\cdot 10^{-2}$ & 45.0  &  5.08$\cdot 10^{-2}$      & \bf 2.0 & 1.14$\cdot 10^{-3}$ \\
		\hline
		\multirow{4}*{\texttt{pref(n,5)}}       & 3200  & 11.4 &    1.08$\cdot 10^{-3}$      & 4.1 & 8.69$\cdot 10^{-4}$ & 45.0 &   1.71$\cdot 10^{-3}$        &\bf  2.0 &  5.44$\cdot 10^{-4}$ \\
		& 6400  & 11.7 &   7.19$\cdot 10^{-3}$        & 3.0 & 2.78$\cdot 10^{-3}$ & 45.0 &  4.84$\cdot 10^{-3}$     & \bf 1.3  & 3.01$\cdot 10^{-4}$ \\
		& 12800 & 12.2 &   1.81$\cdot 10^{-2}$         & 2.5 & 6.44$\cdot 10^{-3}$ & 45.0 &  3.21$\cdot 10^{-2}$        &  \bf 1.1  &  5.57$\cdot 10^{-4}$ \\
		& 25600 & 12.5 &  3.06$\cdot 10^{-2}$  & 2.2 & 1.09$\cdot 10^{-2}$ & 45.0 &  4.80$\cdot 10^{-2}$        &  \bf 1.1  & 1.16$\cdot 10^{-3}$ \\
		\hline
	\end{tabular}
    \end{adjustbox}
\end{table}

The results show that our strategy (iii) requires fewer iterations to converge, and thus smaller CPU times, compared to the other methods considered, especially in the case of large networks. Moreover, the number of iterations required to achieve convergence by methods (i) with $\vvec = \bzero$ and (ii) tends to increase with the size of the network; on the other hand, it decreases with $n$ for method (iii) and for \texttt{pcg($\xvec^{(0)} = \xvec$)}. This behaviour can be explained by observing that, while \texttt{pcg($\xvec^{(0)} = \bzero$)} and series truncation are iterative methods starting from a poor initial guess, namely $\xvec^{(0)}=\bzero$, our approach (iii) and \texttt{pcg($\xvec^{(0)} = \xvec$)} are initialized with $\xvec^{(0)} = \xvec$, which may already be a good approximation to $\xvec^{\cS}$. In fact, thanks to Corollary~\ref{corollary:update_katz}~and~\ref{corollary:update_katz_arc}, we know that if the removed node (resp., edge) has (resp., joins nodes with) a relative small centrality score(s), then $\xvec^{\cS}$ is almost identical to $\xvec$. 
An alternative explanation may also be given, for method (iii) only, in terms of the combinatorics of walks. Indeed, methods (i) with $\vvec = \bzero$ and (iii) compute $\xvec^{\cS}$ by counting from scratch all the walks that \textit{avoid} $s$, which are the ones not effected by the removal; method (iii), on the other hand, focuses on the walks that \textit{visit} $s$. The observed phenomenon may be hence explained because the number of these latter walks relative to network size decreases with an increasing $n$ (except for specific examples such as star graphs) much faster than the number of those that avoid $s$ relative to $n$.

\bigskip
In our second set of experiments, we investigate the accuracy of Algorithm~\ref{alg:update_katz_edge}~and~\ref{alg:update_katz_node}  when a set of nodes ({ resp., }edges) is removed from the graph. We sequentially removed $\lceil n/100 \rceil$ nodes (resp., $\lceil m/100\rceil$ edges) from the synthetic networks in our dataset. After each removal, we compute both the actual value of Katz centrality and its approximation. We perform both a random and a targeted selection of nodes/edges, { which simulate random failures and malicious attacks to the networks, respectively}. The results are averaged over 30 runs.

To assess the effectiveness of the algorithms, we consider the relative 2-norm error $\|\xvec-\widehat{\xvec}\|/\|\xvec\|$ between the actual vector of Katz centralities $\xvec$ and the approximations $\widehat{\xvec}$ computed with Algorithms~\ref{alg:update_katz_edge} and \ref{alg:update_katz_node}. We also look at the \emph{intersection similarity}  of the induced rankings~\cite{isim_article}. The intersection similarity of two ranking vectors $\boldsymbol{\beta}$ and $\boldsymbol{\gamma}$ 
is defined as
\[
	isim_ {p} (\boldsymbol{\beta},\boldsymbol{\gamma}) = \frac{1}{p} \sum_ {i = 1}^ p \frac {| \boldsymbol{\beta}_{1:i} \Delta \boldsymbol{\gamma}_{1:i} |}{2i},
\]
where, using Matlab notation, we denote by $ \boldsymbol{\beta}_{1:i} $ and $ \boldsymbol{\gamma}_{1:i} $ the vectors of length $i$ containing the first $ i $ elements in $\boldsymbol{\beta}$ and $\boldsymbol{\gamma}$, respectively, and we denote by $\Delta $ the symmetric difference operator.  
The intersection similarity is a popular distance often used to compare the rankings given by
two centrality measures; see, e.g., \cite{TC_Benzi_Klymko,multiplex_network_Arrigo_Tudisco,nonlocal_pagerank_Cipolla_Durastante_Tudisco}. If $isim(\boldsymbol{\beta},\boldsymbol{\gamma})\approx 0$, then the two rankings are almost
identical, whereas
$isim(\boldsymbol{\beta},\boldsymbol{\gamma})=1$ represents two completely different rankings.

\begin{figure}[t]
    \includegraphics[width=1.1\textwidth]{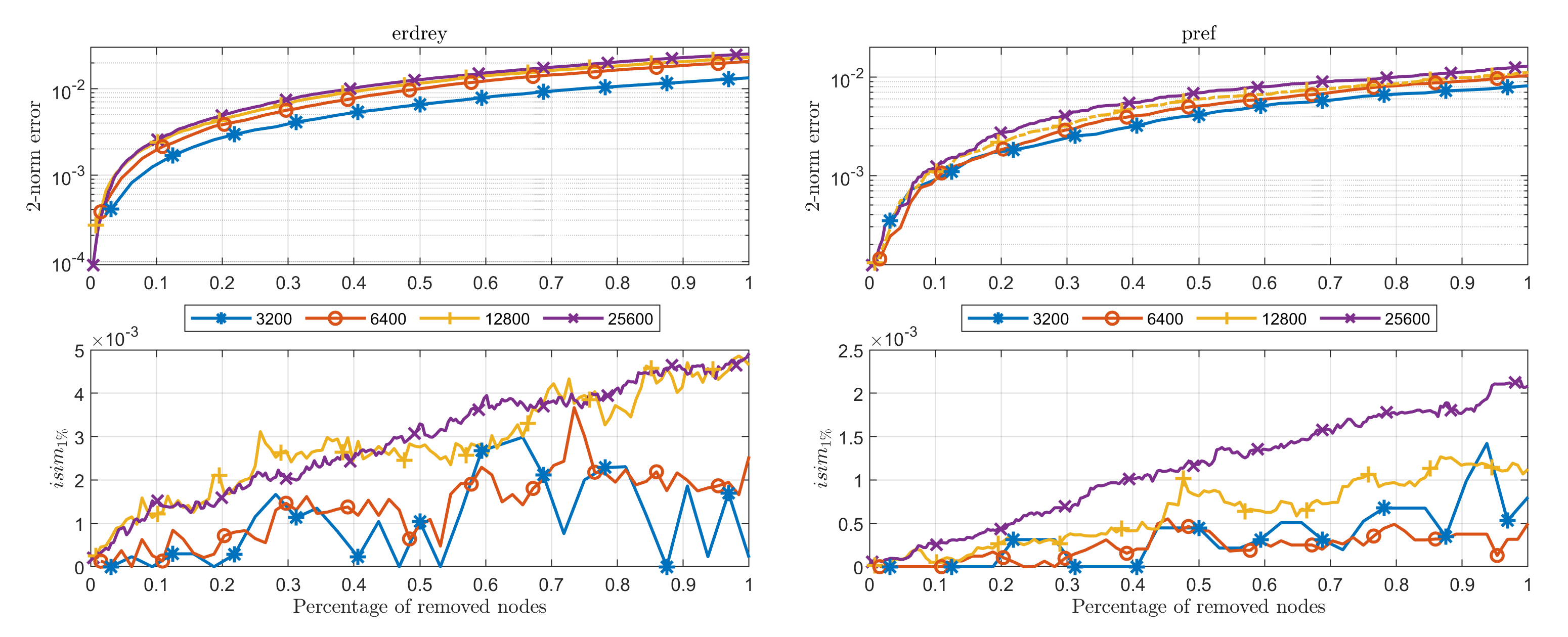}
    \caption{Random removal of a sequence of nodes. Evolution of the relative 2-norm error and $isim_{p}$, where $p=\lceil n/100\rceil$, when sequentially removing nodes from synthetic networks of increasing sizes.}
    \label{fig:error_random_nodes}
\end{figure}

\begin{figure}[t]
    \centering
    \includegraphics[width=1.1\linewidth]{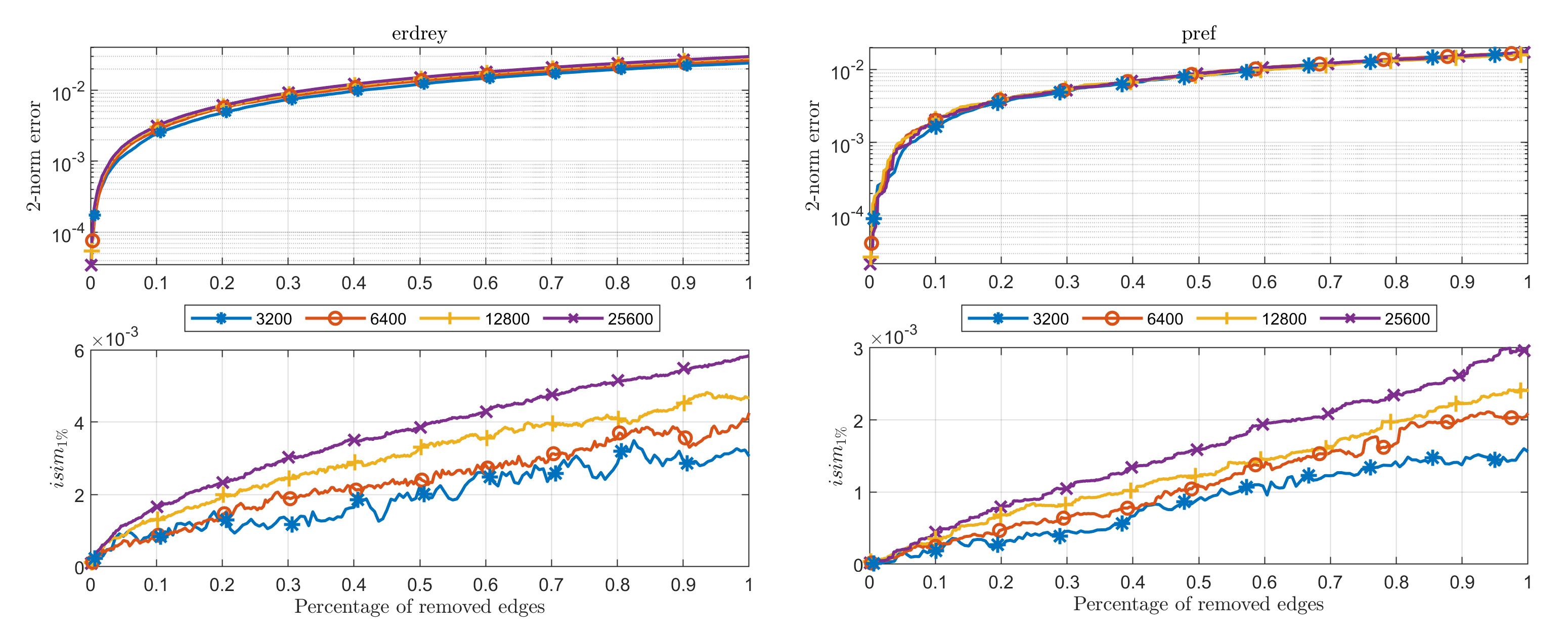}
    \caption{Random removal of a sequence of edges. Evolution of the relative 2-norm error and $isim_{p}$, where $p=\lceil m/100\rceil$, when sequentially removing edges from synthetic networks of increasing sizes.}
\label{fig:error_random_edges}
\end{figure}

Figures~\ref{fig:error_random_nodes}~and~\ref{fig:error_random_edges} plot the evolution of the 2-norm error (top) and of the intersection similarity between the ranking vectors induced by the ground truth $\xvec$ and by the computed approximations $\widehat{\xvec}$ with $p={\lceil n/100\rceil}$ (bottom) as a function of the percentage of nodes/edges being removed. The relative error stabilizes at around $10^{-2}$. The intersection similarity also stabilizes at very low values, indicating that the approximations computed using Algorithm~\ref{alg:update_katz_edge}~and~\ref{alg:update_katz_node} correctly identify the nodes ranked in the top $1\%$ by exact Katz centrality after each removal.

\begin{figure}[t]
    \centering
    \includegraphics[width=1.1\linewidth]{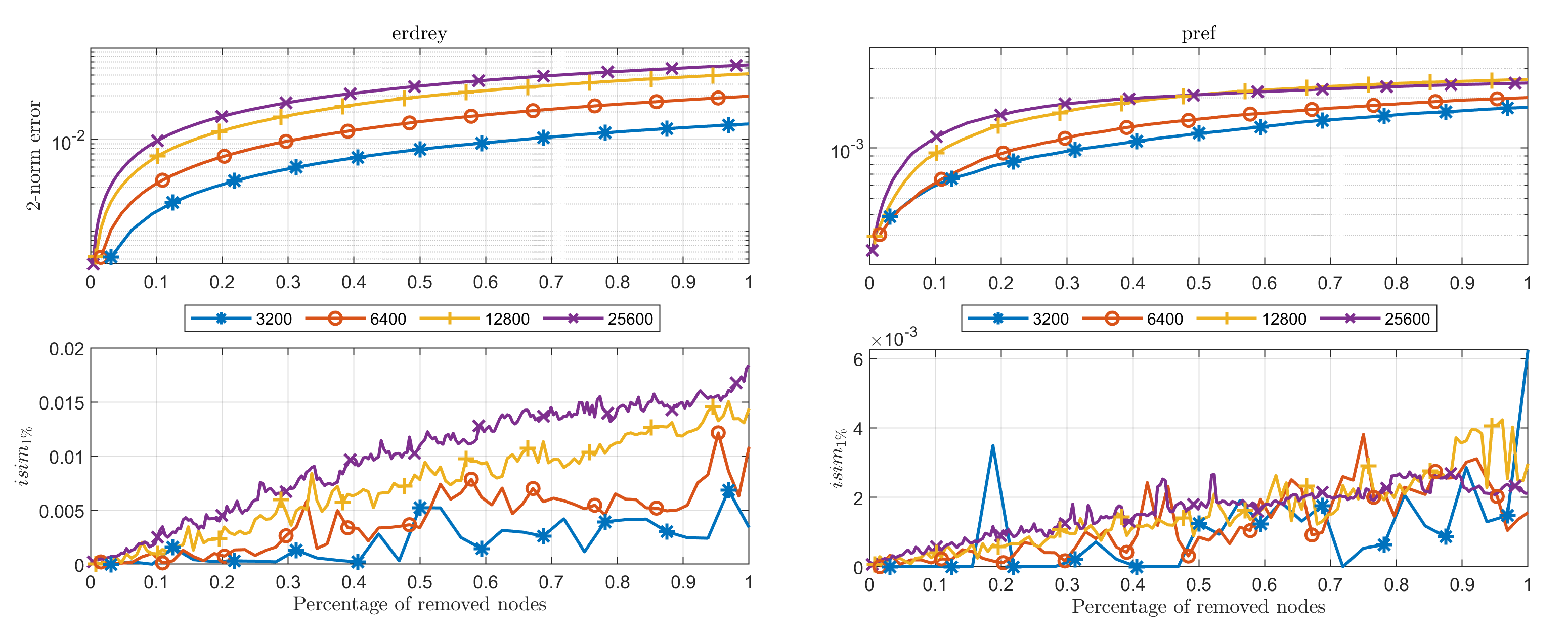}
    \caption{Targeted removal of a sequence of nodes. Evolution of the relative 2-norm error and  $isim_{p}$, where $p=\lceil n/100\rceil$, when sequentially removing $1\%$ of the most important nodes (according to Katz centrality) from synthetic networks of increasing sizes.
    All the values are averaged over $30$ runs.}
    \label{fig:error_most_important_nodes}
\end{figure}

We  now consider the same experiment, but instead of removing elements at random, we remove the top $1\%$ of the highest ranked nodes (according to Katz centrality). 
It is well known that removing these elements will significantly affect the centrality scores and the induced rankings; see, e.g., \cite{On_the_stability_PozzaTudisco}.
We plot the results of our experiment in Figure~\ref{fig:error_most_important_nodes}. These show that our approximations retain a similar error rate even though the centralities changed significantly from their initial value. The values of intersection similarity for these experiments are larger than in the previous one, but still very low. This confirms that our strategies well reproduce both scores and rankings, regardless of node removal strategy.

\bigskip
The last set of experiments examines the behaviour of the bounds in Theorem~\ref{thm:bounds_total_comm} and \ref{thm:bounds_total_comm_arc} when $1\%$ of nodes or edges selected uniformly at random are removed from the two real-world networks in our dataset. { These tests mimic scenarios where a random event shuts down roads/intersections in the Minnesota road network, or isolates a group of routers within the network AS-735 of autonomous systems. As before,} the results are averaged over 30 runs. 

\begin{figure}[t]
    \centering
\includegraphics[width=0.9\linewidth]{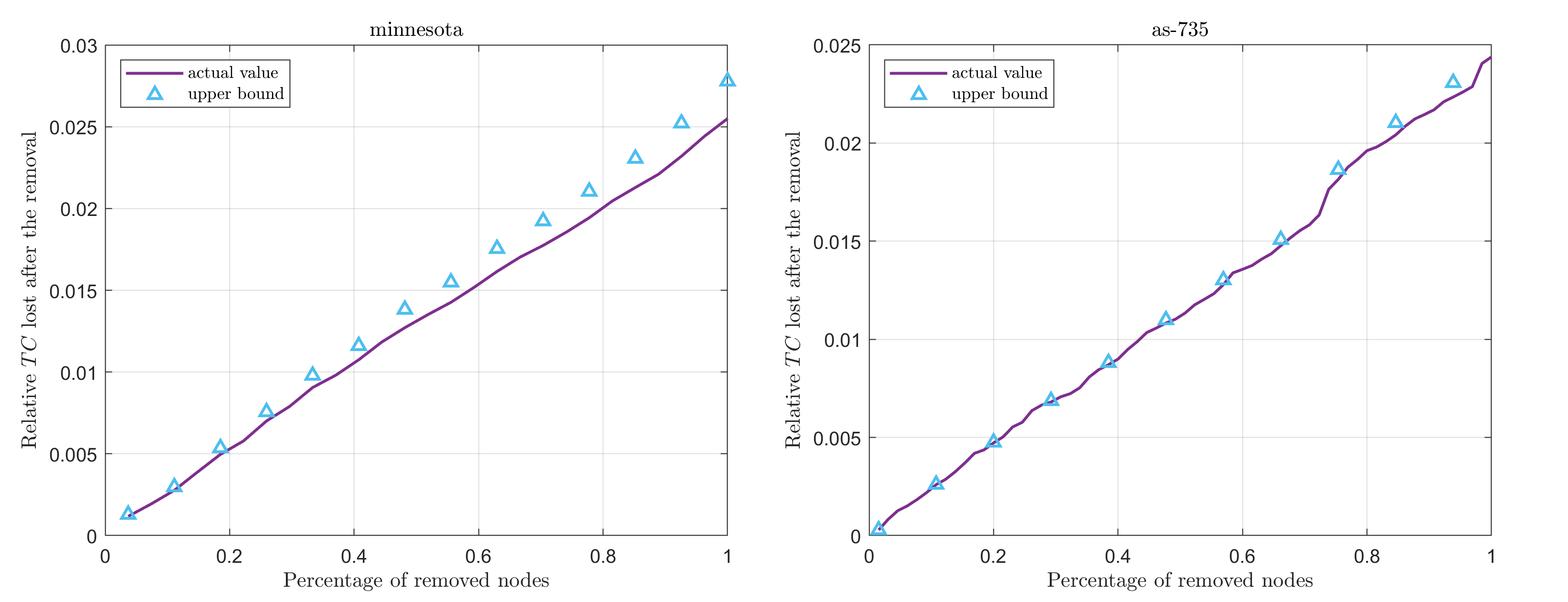}
\caption{Comparison between the bound (\ref{eq:bounds_total_comm})
and the actual value of $TC$ lost (relative to the initial value)
after the removal of $1\%$ of the total nodes form the Minnesota
roads the as-735 connections networks, taken uniformly at random.
The values are averaged over $30$ runs.}
    \label{fig:tc_bounds_node}
\end{figure}

\begin{figure}[t]
    \centering
      \includegraphics[width=0.9\linewidth]{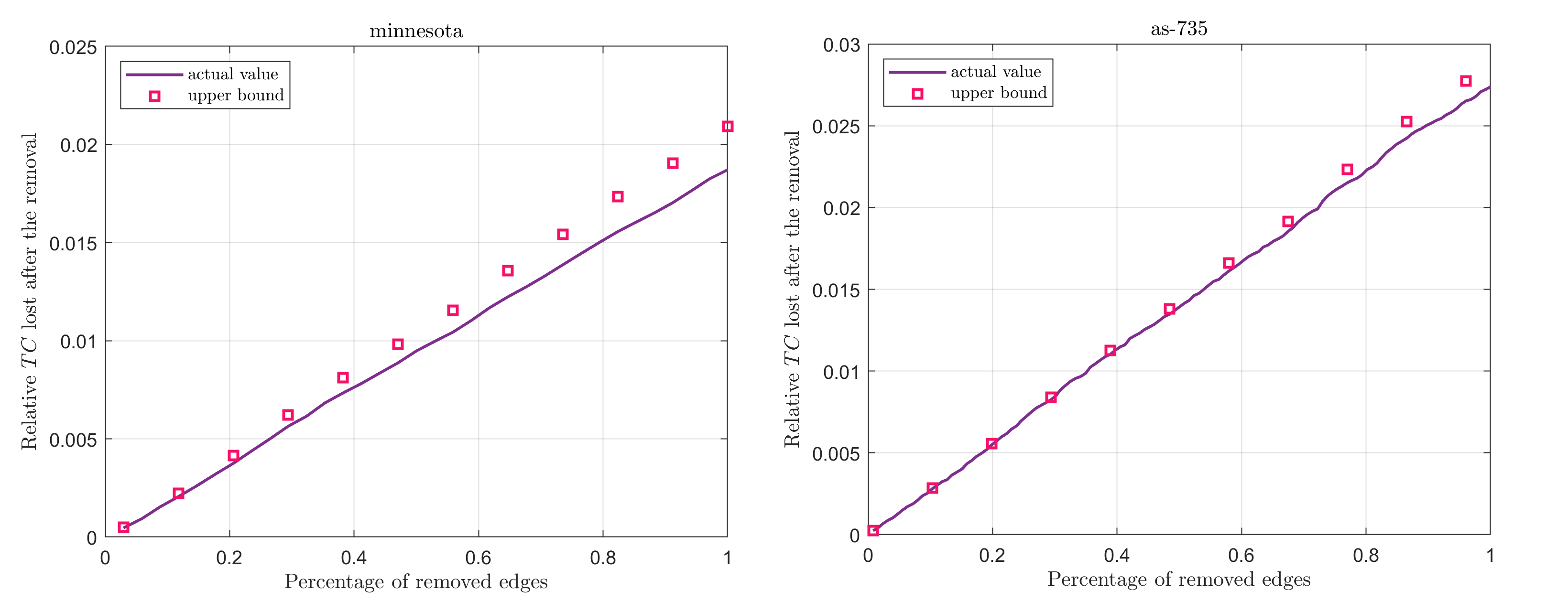}
\caption{Comparison between the bound
\eqref{eq:bounds_total_comm_arc} and the actual value of $TC$ lost
(relative to the initial value) after the removal of $1\%$ of the
total edges form the Minnesota roads and the as-735 connections
networks, taken uniformly at random. The values are averaged over
$30$ runs.}
    \label{fig:tc_bounds_edges}
\end{figure}

The results are displayed in Figures~\ref{fig:tc_bounds_node} and \ref{fig:tc_bounds_edges}, where we plot the change in the total network communicability as nodes (resp., edges) are sequentially removed from the network. On the same plots, we display the upper bounds obtained in \eqref{eq:bounds_total_comm} and \eqref{eq:bounds_total_comm_arc}. The bounds are extremely close to the actual values of total network communicability. A remark is however in order; the sharpness of these bounds strongly depends on the change in the underlying centrality scores. If these change significantly, then using the previous value of Katz centrality to obtain an upper bound will, over time, cause a severe
overestimation of the change in total network communicability. 
This situation may arise when some crucial elements
are removed from the network or when the number of removed elements
becomes large. In these situations, one could mitigate the error by
selectively recomputing the measures after a suitable number of
iterations; see, e.g.,  \cite{ranking_PLOS_Bertaccini_Filippo} for some strategies.

\section{Conclusions and future work}\label{sec:conclusions}
In this paper we derived explicit formulas, expressed in terms of loss of walks,  for the change in Katz centrality when nodes and/or edges are removed from a simple graph. These formulas informed the description of two algorithms to approximate the new centrality vectors after network modification. The theoretical computational cost of the proposed methods is $\mathcal{O}(Lm)$, where $m$ is the number of undirected edges in the original graph and $L$ is a (dynamically chosen) maximum walk-length. We also described upper bounds on the change in total network communicability. 
Numerical results on both synthetic and real-world networks showed that our approximation strategies are more efficient than others available in the literature, including recomputing from scratch and series truncation, and return results that are accurate both in terms of numerics (relative error) and ranking (intersection similarity).

Given these promising results, future work will focus on developing similar strategies for other graph types and centrality measures. 

\section*{Acknowledgements}
The work of DB was partially supported by MIUR Excellence
Department Project MatMod@TOV 2023-2027 awarded to the Department of
Mathematics, University of Rome Tor Vergata and by the Horizon 2020
Project Energy oriented Centre of Excellence: toward exascale for
energy (EoCoE II). DB and AF are members of the Gruppo Nazionale
Calcolo Scientifico-Istituto Nazionale di Alta Matematica
(GNCS-INdAM)"

\bibliographystyle{siam}
\bibliography{biblio}

\end{document}